\documentclass[11pt]{article}
\usepackage[centertags]{amsmath}
\usepackage[normal]{subfigure}
\usepackage{wrapfig}
\usepackage{verbatim}
\usepackage{amsfonts}
\usepackage{amssymb}
\usepackage{amsthm}
\usepackage{amsmath}
\usepackage{newlfont}
\usepackage{graphicx}
\usepackage{epsfig}
\usepackage{amsmath}
\usepackage{color}
\usepackage{multicol}
\usepackage{amssymb,amsmath}
\usepackage{latexsym}
\usepackage[normalem]{ulem}
\usepackage{cancel}

\usepackage{bm}

\usepackage{amsmath, amsthm, amssymb}

\newtheorem{prop}{Proposition}

\DeclareTextFontCommand{\rm}{\rmfamily}

\definecolor{dred}{rgb}{0.92,0,0}
\definecolor{dgreen}{rgb}{0,0.92,0}
\definecolor{dblue}{rgb}{0,0,0.92}
\definecolor{dyellow}{rgb}{0.95,0.95,0}

\usepackage{graphicx}
\usepackage{color}
\usepackage[normal]{subfigure}
\definecolor{dred}{rgb}{0.92,0,0}
\definecolor{dgreen}{rgb}{0,0.92,0}
\definecolor{dblue}{rgb}{0,0,0.92}
\definecolor{dyellow}{rgb}{0.95,0.95,0}

\newtheorem{Rem}{Remark}

\def\ds{\displaystyle}
\def\Bvec{\textbf{B}}
\def\Cvec{\textbf{C}}

\def\Evec{\textbf{E}}
\def\Hvec{\textbf{H}}
\def\Jvec{\textbf{J}}

\def\Yvec{\textbf{Y}}
\def\Xvec{\textbf{X}}
\def\Fvec{\textbf{F}}
\def\Lvec{\textbf{L}}
\def\Svec{\textbf{S}}
\def\nvec{\textbf{n}}
\def\xvec{\mathbf{x}}
\def\yvec{\mathbf{y}}
\def\evec{\textbf{e}}
\def\fvec{\textbf{f}}
\def\uvec{\textbf{u}} 
\def\vvec{\textbf{v}}

\def\wvec{\textbf{w}}

\newcommand{\pa}{\partial}
\newcommand{\curlv}{\textbf{curl}\,}

\newcommand{\gradv}{\textbf{grad}\,} 

\renewcommand{\div}{{\rm{div}}\,}
\graphicspath{/Figures/.}
\graphicspath{./Figures/.}

\def\ds{\displaystyle}

\def\nvec{\textbf{n}}
\def\Z{\mathbb{Z}}
\def\R{\mathbb{R}}

\graphicspath{
{../PhD_Article2_and_Article3/Figures/}
{../Figures/}
}
\title{Numerical solution to the 3D Static Maxwell equations in axisymmetric singular domains with arbitrary data}
\author{
F. Assous\footnote{Ariel University, 40700 Ariel, Israel} ,
I. Raichik \footnote{Bar-Ilan University, 52900 Ramat-Gan, Israel.}
}
\date{}

\begin{document}
\maketitle
\begin{abstract}
We propose a numerical method to solve the three-dimensional static Maxwell equations in a
singular axisymmetric domain, generated by the rotation of a singular polygon around one of its sides.  The mathematical tools and an in-depth study of the problem set in the meridian half-plane are exposed in \cite{AsCLS03}, \cite{CiLa11}. Here, we derive a variational formulation and the corresponding approximation method.  Numerical experiments are proposed, and show that the approach is able to capture the singular part of the solution. This article can also be viewed as a generalization of the Singular Complement Method to three-dimensional axisymmetric problems. 
\end{abstract}

{\bf\noindent keywords}: Maxwell equations, Fourier analysis, Singularities, Axisymmetric geometry, Finite element.

\section{Introduction}\label{intro}

\noindent There is a need to simulate electromagnetic wave phenomena of increasing complexity, leading to the development of  more general and efficient numerical methods. Indeed, a plethora  of engineering problems requires to simulate numerically devices working with or within electromagnetic fields.\\

\noindent  This article is part of the efforts made in the more general  framework of boundary value problems with singularity in their solutions, caused by the presence of {\em geometrical singularities}, that is reentrant corners or edges on the boundary of a domain, or similarly by a change in the type of boundary conditions \cite{Gri85,Gri92}. From a more physical point of view, they are called singularities, since they can generate in their vicinity very strong fields that have to be taken into account, and are very difficult to compute. Moreover, as illustrated in \cite{AsCS00}, the inability to properly handle these singularities may have dramatic consequences on the physical phenomenon one wants to study. \\

\noindent In this context, several authors have proposed to use methods that ``extract" the singular part of the solution near
these singularities,  or to apply mesh refinement toward these singularities, in the case of weak singularities (roughly speaking, that belong to a regular enough space like $H^1$). This allows to construct numerical methods that are able to catch the singular behavior of the solution. The non-matching grid approach is also an interesting alternative \cite{BHS03}.\\

\noindent In this article, we are more specifically concerned with solving three-dimensional Maxwell's equations, that are often used  
to describe the physics of engineering problems, in their static or time-dependent form, sometimes coupled with other equations (see an overview in \cite{AsCL18}).  Moreover, many structures that are to be modeled have a complex three-dimensional geometry that presents a surface with reentrant edges and/or corners, namely singularities.\\

\noindent There exist many methods to solve the Maxwell equations numerically \cite{Monk03}. One can mention the edge finite element method \cite{Nede1,Nede2}, more recently, the class of discontinuous Galerkin method introduced by \cite{HeWa08}, or adaptive finite element method  in two dimensions, as proposed in \cite{BrGS12}. However, it is interesting for some applications to have a continuous approximation of the solutions, for instance when coupling the Maxwell equations with other equations, like Vlasov's one \cite{AsCG00}, \cite{ADHRS93}.\\

\noindent As it is well-known, when solving Maxwell's equations in a non-convex and non-smooth domain with a continuous approximation, the discretized spaces are always included in a closed, strict subspace of the space of real solutions, see the seminal work of \cite{BiSo87,BiSo87a} for theoretical justifications, and more recent developments by Costabel and Dauge (see among others \cite{CoDa00}). Consequently, it is not possible to approximate the singular field and needs special treatment, even for static problems \cite{CGP08}. In this case, mesh refinement techniques fail. The Singular Complement method ({\em SCM}) \cite{AsCLS03, AsCS00} addresses this problem by explicitly adding some singular complements to the space of solutions, see also \cite{AsCG05, ACLL02, AsCS98}.\\

\noindent Numerical solution of three-dimensional boundary value problems in non-convex domains is basically different from the two-dimensional case and is often more difficult. Among many existing methods, Fourier Finite Element Method is an efficient method for solving problems in three-dimensional prismatic or axisymmetric domains, even for other equations, see for instance \cite{BBDH06} for Stokes equations. The method uses the Fourier expansion in one space direction associated to a finite element approach in the other two space dimensions, see, among others \cite{CHQZ88}, \cite{Hein96}, \cite{Gri92}, \cite{MR82}, or \cite{HNW03} for interface problems.\\

\noindent In the present work, that extends the {\em SCM} to three-dimensional axisymmetric singular domains with arbitrary data, 
we consider a situation in which the three-dimensional (3D) Maxwell equations can be reformulated as two-dimensional (2D) models. This principle was also derived in \cite{CJKLZ05,CJKLZ05I} for the Poisson equation in a prismatic or axisymmetric geometry.
More precisely, the computational domain boils down to a subset of $\R^2$, with respect to the cylindrical system of coordinates. Nevertheless, the electric and magnetic fields, and other vector quantities, still belong to $\R^3$. Hence, the electromagnetic field is the solution to an infinite set of 2D equations, and as a result a set of 2D variational formulations, obtained by Fourier analysis.\\

\noindent This paper is organized as follows: in a first section, we recall the Maxwell equations and their formulation in an axisymmetric domain. Then we present the principle of the 2D space reduction, based on the use of a Fourier transform in $\theta$. This reduces 3D Maxwell's equations to a series of 2D Maxwell's equations, depending on the Fourier variable $k$. This allows us to compute the 3D solution by solving several 2D problems, each one depending on $k$. Even if the solution remains singular for each $k$ in the 2D domain, we will be able to  decompose it into a regular and a singular part (see Section \ref{DecompRegSin}). The regular part belongs to a regular space and will be computed  by a standard finite element method. The singular part, that belongs to a  finite-dimensional subspace, will be handled following the same principle as in the {\em SCM}. This is the subject of Section \ref{ComputSingu}. In the last Section, numerical examples are proposed to illustrate the feasibility of the method.\\

\noindent  In the remainder of this paper, we write vector fields or spaces with boldface. Similarly, names of function spaces of 
scalar fields usually begin by an italic letter, whereas they begin by a bold letter for spaces of vector fields (for 
instance, $\Lvec^2(\Omega)=L^2(\Omega)^{3}$ or $L^2(\Omega)^{2}$).
\section{Maxwell's equations in an axisymmetric domain}\label{MaxEquat}
\subsection{The static Maxwell equations}
\noindent Let $\Omega$ be a bounded and simply connected Lipschitz domain in $\mathbb{R}^{3}$, $\Gamma$ its boundary, assumed for simplicity to be a connected boundary, and $\nvec$ the unit outward normal to $\Gamma$. If we let $c$ and $\varepsilon_0$ be, respectively,  the speed of light and the dielectric permittivity, the time-dependent Maxwell equations in vacuum read in $\Omega$,
\begin{eqnarray}
&&\frac{\partial { \Evec}}{\partial t}-c^2{{ \curlv}\,}{ \Bvec}=-\frac{1}{\varepsilon_0}{\Jvec}, \label{ampeq3D}\\ 
&&\frac{\partial { \Bvec}}{\partial t}+{{ \curlv}\,}{\Evec}=0,\label{faradeq3D}\\
&&{{\div}\,}{ \Evec}=\frac{{\rho}}{\varepsilon_0},\label{couleq3D}\\ 
&&\div\,\Bvec=0\,,\label{monopeq3D}
\end{eqnarray}
where $ \Evec$ is the electric field,  $\Bvec$ the magnetic flux density, $\rho$ and $\Jvec$ the charge and current
densities. These quantities depend on the space variable $\xvec$ and on the time variable $t$.\\

\noindent These equations are supplemented with appropriate boundary conditions. In this article, we assume that the boundary $\Gamma$ is a perfect conductor, so that the electromagnetic field satisfies
\begin{equation}
\label{BC3D}
\Evec\times \nvec=0 \quad \mbox{ and }\quad  \Bvec\cdot \nvec=0 \quad \mbox{ on the boundary  } \Gamma\,.
\end{equation}

\noindent Since we are interested in the static Maxwell equations, we consider problems and solutions that are time-independent, namely static equations.  In other words, we assume that the explicit time-dependence $\partial / \partial t$ of the electromagnetic field in Maxwell's equations vanishes. With non-vanishing charge and current densities, this assumption yields there are two div-curl problems, depending on the boundary condition.\\

\noindent  The first one is, for a mean zero value right-hand side $\fvec_\Evec$ in $\Lvec^2(\Omega)$, such that $\div\fvec_\Evec=0$ and $\fvec_\Evec\cdot\nvec_{|\Gamma}=0$, and for a right-hand side $g_\Evec$ in $L^2(\Omega)$:\\
{\em Find $\Evec\in \Lvec^2(\Omega)$ such that}
 \begin{eqnarray}
&&\curlv\Evec = \fvec_\Evec\mbox{ in }\Omega,\label{es1}\\
&&\div\Evec = g_\Evec\mbox{ in }\Omega,\label{es2}\\
&&\Evec\times\nvec_{|\Gamma} = 0.\label{es3}
 \end{eqnarray}
 The boundary condition on $\fvec_\Evec$ is imposed by the condition (\ref{es3}) (cf. \cite{GiRa86}). In order to prove the existence and uniqueness of the solution $\Evec$ to (\ref{es1})-(\ref{es3}), a possible way is to reformulate these equations as a saddle-point formulation, and to check that the Lagrange multiplier is equal to 0 (see \cite{CiZo97}, \cite{AsCL18} Chap.6 for details). Assuming the connectivity of the boundary $\Gamma$ is required here, since the use of the saddle point approach needs to use a Friedrichs-type inequality of the form $\|\vvec\|_0 \leq C \| \curlv \vvec\|_0$. Equivalently, Eqs. (\ref{es1})-(\ref{es3}) can represent the stationary problem associated with Maxwell's equations, namely the quasi-electrostatic problem. This amounts to assuming that the time-dependent parts are known, and that $\fvec_\Evec=\ds\frac{\partial { \Bvec}}{\partial t}$ and $g_\Evec=\ds\frac{{\rho}}{\varepsilon_0}$.\\

\noindent With analogous notations, the second div-curl problem is, for a given $\fvec_\Bvec$ in $\Lvec^2(\Omega)$ such that $\div\fvec_\Bvec=0$, and for a mean zero value $g_\Bvec$ in $L^2(\Omega)$:\\
 {\em Find $\Bvec \in \Lvec^2(\Omega)$ such that}
 \begin{eqnarray}
&&\curlv\Bvec = \fvec_\Bvec\mbox{ in }\Omega,\label{bs1}\\
&&\div\Bvec = g_\Bvec\mbox{ in }\Omega,\label{bs2}\\
&&\Bvec\cdot\nvec_{|\Gamma} = 0.\label{bs3}
 \end{eqnarray}
Similarly, this can model the quasi-magnetostatic Maxwell's equations by assuming $\fvec_\Bvec=-\ds\frac{1}{c^2}\frac{\partial { \Evec}}{\partial t} + \frac{1}{c^2\,\varepsilon_0}{\Jvec}$ and $g_\Bvec=0$. The fact that $g_\Bvec$ has a mean zero value stems from (\ref{bs3}). The existence and uniqueness of $\Bvec$ can also be inferred by using a saddle-point approach. In both cases, the existence and uniqueness result can be achieved thanks to the Weber 
inequality \cite{Webe80}, see details in \cite{AsCL02} or in \cite{AsCL18} Chap.6.

\subsection{Formulation in an axisymmetric domain}
\begin{figure}[htbp!]
\begin{tabular}{lr}
{
\includegraphics[scale=0.25]{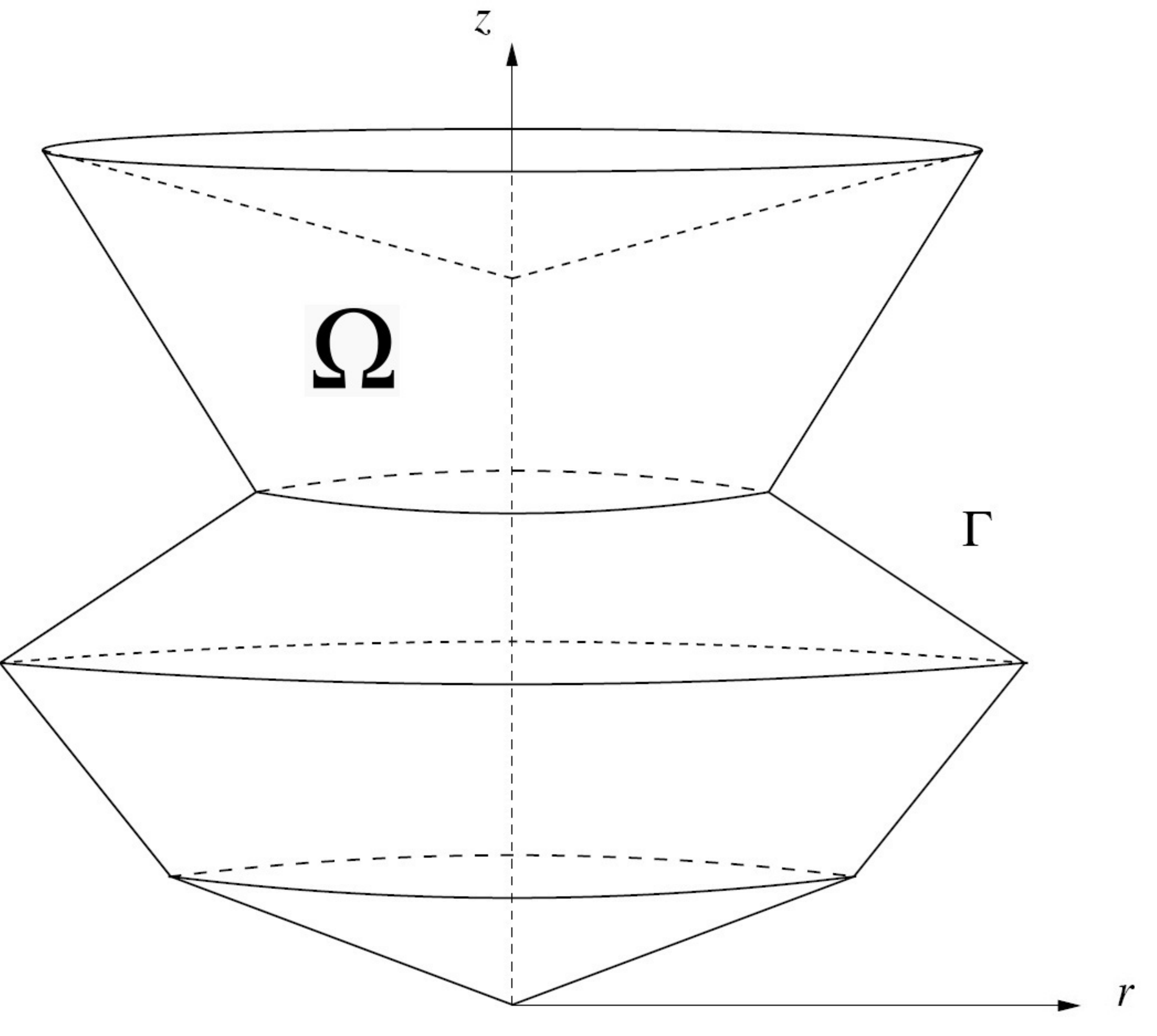}
}
&
\hspace*{2.cm}
{
\includegraphics[scale=0.5]{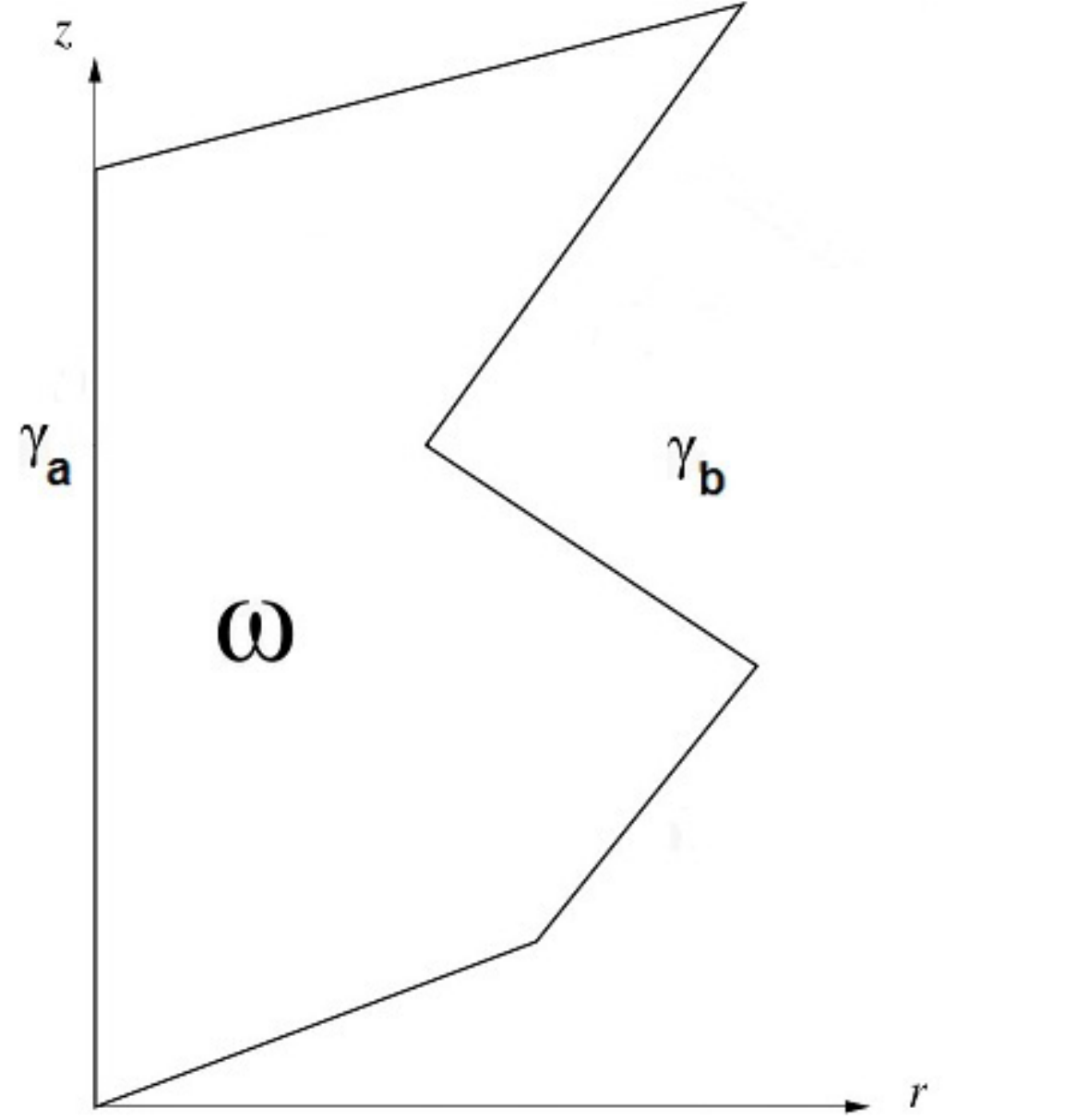}
}
\end{tabular}
 \caption{The $\Omega$ and $\omega$ domains.}
  \label{fig:domain}
\end{figure}

Now we make the supplementary assumption that $\Omega$ is an axisymmetric domain, limited by the surface of revolution $\Gamma$. We denote by $\omega$ and $\gamma _b$  their intersections with a meridian half-plane (see Figure \ref{fig:domain}). One has $\partial \omega :=\gamma =\gamma _a \cup \gamma _b$, where either $\gamma _a= \emptyset$ when $\gamma_b$  is a closed contour (i.e. $\Omega$ does not contain the axis), or $\gamma_a$  is the segment of the axis lying between the extremities of $\gamma_b$. 
The natural coordinates for this domain are the cylindrical coordinates $(r,\theta,z)$, with the basis vectors $(\evec_r, \evec_\theta,\evec_z)$. A meridian half-plane is defined by the equation $\theta=$constant,  and $(r,z)$ are Cartesian coordinates in this half-plane.\\

\noindent However, although the domain $\Omega$ is assumed to be axisymmetric,  the symmetry of revolution {\em is not assumed} for the data. In these conditions, the problem can not be reduced to a two-dimensional one by assuming that $\partial/\partial \theta=0$, as made for instance in \cite{AsCLS03}. We continue here to deal with a three-dimensional problem.

\noindent In these conditions, one can obtain the expressions of the static Maxwell equations simply by replacing into (\ref{es1}-\ref{es3}) and (\ref{bs1}-\ref{bs3})  
the operators $\div$ and $\curlv$ by their cylindrical counterparts  in the cylindrical coordinates $(r,\theta,z)$, with the basis vectors $(\evec_r, \evec_\theta,\evec_z)$, defined by
\begin{eqnarray}
&&\div \uvec= \frac{1}{r}\frac{\partial }{\partial r}\left ( ru_r \right )+\frac{1}{r}\frac{\partial u_\theta}{\partial \theta}+\frac{\partial u_z}{\partial z} \label{div3d}\\
&&\curlv \uvec=\left (\frac{1}{r}\frac{\partial u_z }{\partial \theta }-\frac{\partial u_\theta }{\partial z }   \right )\evec_r+\left ( \frac{\partial u_r }{\partial z } -\frac{\partial u_z }{\partial r }\right )\evec_\theta+\frac{1}{r}\left ( \frac{\partial }{\partial r}\left ( ru_\theta \right )-\frac{\partial u_r }{\partial \theta } \right )\evec_z \label{curl3d}
\end{eqnarray}
Similarly, the gradient operator in cylindrical coordinates is defined by 
\begin{eqnarray}
&&\gradv f=\frac{\partial f}{\partial r}\evec_r+\frac{1}{r}\frac{\partial f }{\partial \theta }\evec_\theta+\frac{\partial f }{\partial z }\evec_z \label{grad3d}
\end{eqnarray}

\subsection{Variational formulations in 3D}
We now introduce the variational formulations of the problem, which can be applied independently of the (non) convexity of the domain $\Omega$. Let us define the function spaces, with classical notations: the usual norm and scalar product of 
$\Lvec^2(\Omega)$ are denoted by $\|\cdot\|_0$ and $(\cdot,\cdot)$ respectively. We shall also need to use the following Sobolev spaces and norms
\begin{eqnarray*}
&&\Hvec(\curlv,\Omega)=\{\vvec \in \Lvec^2(\Omega),\, \curlv\vvec \in 
\Lvec^2(\Omega)\}\,,\qquad\|\vvec\|^2_{\curlv}=\|\vvec\|^2_{0}+\|\curlv\vvec\|^2_{0}\,,\\
&&\Hvec(\div,\Omega)=\{\vvec \in \Lvec^2(\Omega),\, \div\vvec \in 
L^2(\Omega)  
\}\,,\qquad\quad\,\|\vvec\|^2_{\div}=\|\vvec\|^2_{0}+\|\div\vvec\|^2_{0}\,,\\
&&\Hvec^1(\Omega)=\{\vvec \in \Lvec^2(\Omega), \,\gradv\vvec \in 
\Lvec^2(\Omega)\}\,,\!\qquad\qquad\|\vvec\|^2_{1}=\|\vvec\|^2_{0}+\|\gradv\vvec\|^2_{0}\,.
\end{eqnarray*}
We introduce likewise
$$
\Hvec_0\left ( \curlv, \Omega\right )=\left \{ \vvec \in \Hvec\left ( \curlv, \Omega\right ) : \vvec \times \nvec|_{\Gamma}=0 \right \} 
$$
and
$$
\Hvec_0\left ( \div, \Omega\right )=\left \{ \vvec \in \Hvec \left( \div, \Omega\right ) : \vvec \cdot \nvec|_{\Gamma}=0 \right \} \,.
$$
The electric and magnetic field naturally belongs respectively to the spaces 
$$
\Xvec\left(\Omega \right ) =\Hvec_0\left ( \curlv, \Omega\right )\cap\Hvec\left ( \div, \Omega\right ) 
\,\,\mbox{ and } \,\,
\Yvec\left(\Omega \right ) =\Hvec\left ( \curlv, \Omega\right )\cap\Hvec_0\left ( \div, \Omega\right ) \,.
$$
The spaces $\Xvec\left ( \Omega \right )$ and $\Yvec\left ( \Omega \right ) $ are compactly embedded in $\Lvec^2(\Omega)$ \cite{Webe80}, \cite{Cost90}. Consequently, when the boundary $\Gamma$ is connected, one can define an equivalent scalar product and norm on $\Xvec\left ( \Omega \right )$ and $\Yvec\left ( \Omega \right )$ as 
$$
a\left (\uvec,\vvec \right ):=\left (\curlv\uvec, \curlv\vvec \right) +\left (\div\uvec, \div\vvec \right), \quad \|\uvec\|_{\Xvec}= \|\uvec\|_{\Yvec}:=a\left (\uvec,\uvec \right )^{1/2}
$$
In other words, the $\Lvec^2$-norm is uniformly bounded by the $\Xvec$ and the $\Yvec$ norm for elements of $\Xvec\left ( \Omega \right )$ and $\Yvec\left ( \Omega \right )$ respectively. This is the Weber inequality, that basically claims that in $\Xvec\left ( \Omega \right )$ or in $\Yvec\left ( \Omega \right )$,  the semi-norm $\uvec \longrightarrow (\| \curlv \uvec\|_0^2 +\|\div\uvec\|_0^2)^{1/2}$ is a norm equivalent to the canonical one.\\

\noindent We have now to derive the (augmented) variational formulations associated to these problems. Following a classical approach, we first take the dot product of equations (\ref{es1}) (resp: (\ref{bs1})) by $\curlv \Fvec$, $\Fvec \in \Xvec \left(\Omega \right)$ (resp: $\curlv \Cvec$, $\Cvec \in \Yvec \left(\Omega \right)$) and integrate over $\Omega$, then add the variational form of the divergence equation for $\Evec$ (resp: $\Bvec$). This gives the variational formulations:\\
 {\em  Find $\Evec\in \Xvec \left ( \Omega \right )$ such that:}
\begin{equation}
\label{varE3D}
  a\left ( \Evec,\Fvec \right )=\left (\fvec_\Evec,\curlv \Fvec \right )+\left (g_\Evec,\div \Fvec \right ), \forall \Fvec\in \Xvec \left(\Omega \right)\,,
\end{equation}
and similarly, for the magnetic field,\\
{\em  Find $\Bvec\in \Yvec \left ( \Omega \right )$ such that:}
\begin{equation}
\label{varB3D}
a\left ( \Bvec,\Cvec \right )=\left (\fvec_\Bvec,\curlv \Cvec \right )+\left (g_\Bvec,\div \Cvec \right ), \forall \Cvec\in \Yvec \left(\Omega \right)\,.
\end{equation}

\noindent Existence, uniqueness and continuous dependence with respect to the data of these variational formulations follow from the application of usual techniques, see for instance \cite{AsCL18}.

\section{Principle of two-dimensional space reduction}\label{DimReduc}
Since we consider non axisymmetric data, we can not perform $\pa/\pa \theta=0$ to reduce the 3D space problem to a 2D space one. However, we will use that the domain $\Omega$ is axisymmetric. The scalar and vector fields defined on $\Omega$ will be characterized through their Fourier series in $\theta$, the coefficients of which are functions defined on $\omega$. Note that such a technique together with the Fourier-Laplace transform is also used for stability analysis of numerical schemes \cite{ChMo10} solving Maxwell's equations. 
Let us also emphasize that the time dependent part of the problem is not involved here.  What is explained below is the {\em principle} of the two-dimensional space reduction. For this reason, we do not mention the time variable in the Fourier series, which can be easily added the case occurring. For instance, we will consider for a given function $w(r,\theta,z)$ (resp: for a vector field $\wvec(r,\theta,z)$), the Fourier expansion
$$
w\left ( r,\theta,z \right )=\frac{1}{\sqrt{2\pi}}\sum_{k\in\mathbb{Z}}w^{k}\left (r,z  \right )e^{ik\theta} ,$$ resp.$$ \wvec\left ( r,\theta,z \right )=\frac{1}{\sqrt{2\pi}}\sum_{k\in\mathbb{Z}}\wvec^{k}\left (r,z  \right )e^{ik\theta} $$

\noindent and the truncated Fourier expansion of $\wvec$ at order $N$

\begin{equation}
\label{FourierTrunc}
 \wvec^{\left [ N \right ]}\left ( r,\theta,z \right )=\frac{1}{\sqrt{2\pi}}\sum_{k=-N}^{N}\wvec^{k}\left (r,z  \right )e^{ik\theta}.
\end{equation}

\noindent We also consider the weighted Lebesgue space

$$L_{1}^{2}(\omega):=\left \{ w\mbox{  measurable on }  \omega:\iint_{\omega}\left | w(r,z) \right |^{2}rdrdz< \infty  \right \}$$
which is the space of Fourier coefficients (at all modes) of functions in $L^{2} \left ( \Omega \right ) $.

\noindent Let us now examine the space of relevant Fourier coefficients for the electromagnetic fields. One easily checks that for $w \in H^{1}\left ( \Omega \right )$, resp. $w\in L^{2}\left ( \Omega \right )$ such that 
$\Delta w \in L^{2}\left ( \Omega) \right.$, there holds:
$$
\gradv w=\frac{1}{\sqrt{2\pi}}\sum_{k\in\mathbb{Z}}\gradv_kw^{k}e^{ik\theta}, \mbox{ resp. }  \Delta w=\frac{1}{\sqrt{2\pi}}\sum_{k\in\mathbb{Z}}\Delta_kw^{k}e^{ik\theta}\,,
$$ 

\noindent  while for $\wvec \in \Hvec \left( \div;\Omega \right)$, resp. $\Hvec \left( \curlv;\Omega \right)$:

$$ 
\div \wvec=\frac{1}{\sqrt{2\pi}}\sum_{k\in\mathbb{Z}}\div_k\wvec^{k}e^{ik\theta}\, \mbox{ resp.  } \curlv \wvec=\frac{1}{\sqrt{2\pi}}\sum_{k\in\mathbb{Z}}\curlv_k\wvec^{k}e^{ik\theta} .
$$

\noindent Above, the operators for the mode $k$ are defined as:
\begin{equation}
\label{eq-oper-k}
\begin{matrix}
\gradv_k:=\ds\frac{\partial w}{\partial r}\evec_r+\frac{ik}{r}w\evec_\theta+\frac{\partial w}{\partial z}\evec_z ;& \Delta_k w:=\ds\frac{1}{r}\frac{\partial }{\partial r}\left ( r\frac{\partial w}{\partial r}\right ) -\frac{k^{2}}{r^{2}}w+\frac{\partial^2 w}{\partial z^2}; \\ 
 \div_k \wvec:=\frac{1}{r}\frac{\partial \left ( rw_r \right )}{\partial r}+\frac{ik}{r}w_\theta+\frac{\partial w_z}{\partial z} ;& \left ( \curlv_k\wvec \right )_{r}:=\frac{ik}{r}w_z-\frac{\partial w_\theta}{\partial z} ;\\ 
\left ( \curlv_k \wvec \right )_{\theta}:=\frac{\partial w_r}{\partial z}-\frac{\partial w_z}{\partial r};& \left ( \curlv_k \wvec \right )_{z}:=\frac{1}{r}\left ( \frac{\partial \left ( rw_\theta \right )}{\partial r} -ikw_r\right ).
\end{matrix}
\end{equation}
As explained in \cite{BeDM99}, the regularity of the function $w$ and of the components of $\wvec$ in the {\em ad hoc} Sobolev spaces are characterized by the regularity  of the Fourier components $w^k$ and $\wvec^k$, for $k \in \mathbb{Z}$.\\

\noindent As a consequence,  a function $\vvec$ belongs to $\Xvec(\Omega)$ if and only if, for all $k \in \mathbb{Z}$, its Fourier coefficients 
$\vvec^k$ belong to the space $\Xvec_{(k)}(\omega)$ defined by
$$
\Xvec_{(k)}(\omega)= \{\vvec^k \in \Lvec_1^2(\omega),\: \curlv_k \vvec^k \in \Lvec_1^2(\omega)\,, \div_k \vvec^k \in L_1^2(\omega)\,,\vvec^k \times \nvec_\mid \gamma_b=0\}
$$
with  $\ds\sum_{k \in \Z} \|\vvec^k \|_{\Xvec_{(k)}(\omega)}^2 < \infty$.\\

\noindent Similarly, we introduce the analogous of $\Yvec(\Omega)$ for the Fourier coefficients, namely 
$$
\Yvec_{(k)}(\omega)= \{\vvec^k \in \Lvec_1^2(\omega),\: \curlv_k \vvec^k \in \Lvec_1^2(\omega)\,, \div_k \vvec^k \in L_1^2(\omega)\,,\vvec^k \cdot \nvec_\mid \gamma_b=0\}
$$
and  $\ds\sum_{k \in \Z} \|\vvec^k \|_{\Xvec_{(k)}(\omega)}^2 < \infty$.\\

\noindent An important property concerning these spaces is proved in \cite{CiLa11}, Prop.2.9:
\begin{prop}
\label{propstab}
The spaces $\Xvec_{(k)}(\omega)$ and $\Yvec_{(k)}(\omega)$  are independent of $k$, for $|k| \geq 2.$
\end{prop}
\noindent As a consequence, it will be sufficient to compute the singular subspaces only for the modes $|k| \leq 2$, while the modes $\pm 2$ will be used to compute all the higher modes $|k| > 2$.

\subsection{Variational formulation in 2D  for each $k$}
Our aim is now to apply this space dimension reduction to the 3D equations, to derive the corresponding 2D formulations satisfied by Fourier coefficients ($\Evec_k$, $\Bvec_k$) for each mode $k$.\\

\noindent More precisely,  we use the linearity of Maxwell's equations (\ref{es1}-\ref{es3}) and  (\ref{bs1}-\ref{bs3})  (or equivalently of their variational formulations) together with the orthogonality
of the different Fourier modes in $\Lvec^2(\omega)$. This implies that the Fourier coefficients $\Evec_k$ and $\Bvec_k$
of $\Evec$ and $\Bvec$ are solutions to variational formulations similar to (\ref{varE3D}) and (\ref{varB3D}), with the operators $ \curlv_k$ and $\div_k$.  
Consequently, introducing the operator $a_k(\cdot,\cdot)$ as follows
\begin{equation}
\label{akexpression}
\begin{array}{l}
 \ds a_{k}(\uvec,\vvec)=
\left (\curlv_k \uvec  ,\curlv_k \vvec  \right ) + \left (\div_k \uvec,\div_k \vvec  \right ),\\ 
\end{array}
\end{equation}
we get that each mode $ \Evec^{k}$ is the solution to the problem:\\
 {\em find $ \Evec^{k} \in \Xvec_{(k)}(\omega)$ such that, for all  $\Fvec \in \Xvec_{(k)}(\omega)$  :} 
 \begin{equation}
\label{akexpressionE}
a_{k}\left (\Evec^{k},\Fvec  \right )= \left (\fvec^k_\Evec,\curlv_k \Fvec \right )+\left (g^k_\Evec,\div_k \Fvec \right )\,,
\end{equation}
where $\fvec^k_\Evec$ and $g^k_\Evec$ denote the Fourier coefficients of the right-hand sides $\fvec_\Evec$ and $g_\Evec$ respectively, that depend only on $(r,z)$.\\

\noindent Similarly,  the magnetic field $\Bvec$ being solution to (\ref{varB3D}), its Fourier coefficients $\Bvec_k$ satisfy the formulation,  for each mode $k$:\\
 {\em find $ \Bvec^{k} \in \Yvec_{(k)}(\omega)$ such that, for all  $\Cvec \in \Yvec_{(k)}(\omega)$  :} 
\begin{equation}
\label{akexpressionB}
a_{k}\left (\Bvec^{k},\Cvec  \right )= \left (\fvec^k_\Bvec,\curlv_k \Cvec \right )+\left (g^k_\Bvec,\div_k \Cvec \right ) \,.
\end{equation}
Here again,  $\fvec^k_\Bvec$ and $g^k_\Bvec$ denote the Fourier coefficients of the right-hand sides involved in the equations of the magnetic field.\\

\noindent For an analysis of the truncation error of the Fourier expansion, we refer the interested reader to  \cite{CiLa11}, \cite{Nkem07}. Basically, the convergence of the truncated solution (see (\ref{FourierTrunc})) $\Evec^{[N]}$  toward $\Evec$ is in $N^{-2s}(\|\fvec_\Evec\|^2 + \|g_\Evec\|^2 )$ for a given norm, the value of $s>1/2$ depending on the regularity of the right-hand sides $\fvec_\Evec, g_\Evec$. Similar results are available for $\Bvec^{[N]}$.
\section{Decomposition in regular/singular parts}\label{DecompRegSin}
Up to now, in the same spirit as in \cite{MR82}, we have reduced the 3D Maxwell equations to a series of 2D Maxwell equations, depending on the Fourier variable $k$.  Nevertheless, the two dimensional domain $\omega$ being singular (see Figure \ref{fig:domain}), we have now to deal with this singularity. The construction of the numerical method will be based on theoretical results proved in \cite{AsCL02}, for the case $k=0$, and in \cite{CiLa11} for the general case.\\

\noindent  For our purpose, we first consider, for each Fourier mode $k$,  the weighted Sobolev space  $\Hvec_{(k)}^{1}(\omega)$ that contains functions $\vvec^k \in \Lvec_1^2(\omega)$ such that  $\gradv_k \in  \Lvec_1^2(\omega)$. Then,  we introduce the regularized spaces $\Xvec^R_{(k)}(\omega)$ and $\Yvec^R_{(k)}(\omega)$ subspaces of $\Hvec_{(k)}^1(\omega)$, defined by:
$$
 \Xvec^R_{(k)}(\omega):=\Xvec_{(k)}(\omega)\cap \Hvec_{(k)}^{1}(\omega) \quad \mbox{ and } \quad \Yvec^R_{(k)}(\omega):=\Yvec_{(k)}(\omega)\cap \Hvec_{(k)}^{1}(\omega)\,.
$$

\noindent We then have the following property (\cite{CiLa11} Lemma 6.2 and \S 6)
\begin{prop}
The regularized spaces $\Xvec^R_{(k)}(\omega)$ and $\Yvec^R_{(k)}(\omega)$ are closed, respectively, within $\Xvec_{(k)}(\omega)$ and $\Yvec_{(k)}(\omega)$.
\end{prop}
\noindent In these conditions, for a singular domain, the spaces of solution $\Xvec_{(k)}(\omega)$ and $\Yvec_{(k)}(\omega)$ can be decomposed in 
$$
\Xvec_{(k)}(\omega)=\Xvec^R_{(k)}(\omega)\oplus \Xvec^S_{(k)}(\omega) \qquad \mbox{ and } \qquad \Yvec_{(k)}(\omega)=\Yvec^R_{(k)}(\omega)\oplus \Yvec^S_{(k)}(\omega)\,.
$$
The subspaces $\Xvec^R_{(k)}(\omega)$ and $\Yvec^R_{(k)}(\omega)$ are the spaces of solutions in case of a regular domain, whereas 
$\Xvec^S_{(k)}(\omega)$  and $\Yvec^S_{(k)}(\omega)$ are singular subspaces, equal to $\{ 0 \}$ for a regular domain. As a consequence, the electromagnetic field solution
 $\Evec^{k}$ and $\Bvec^{k}$ can also be decomposed into a regular  and a singular part, says 
\begin{equation}
\label{decompEBk}
\Evec^{k} = \Evec_R^{k} + \Evec_S^{k}, \qquad 
\Bvec^{k} = \Bvec_R^{k} + \Bvec_S^{k} \,.
\end{equation}
Moreover, these singular subspaces are of finite dimension, the dimension of which depending on the number $N_S$ of singularities in the domain $\omega$.\\

\noindent Hence, one can compute a numerical approximation of $\Evec_R^{k}$ and $\Bvec_R^{k}$ by a standard numerical method, for instance a $P_1$-conforming finite element method. The difficulty coming from the {\em singular} parts $\Evec_S^{k}$  and $\Bvec_S^{k}$, we have now to derive a way to characterize these singular fields. For this purpose, we refer to the following property, see for instance \cite{AsCL03}, Prop. 3.2 of \cite{AsCLS03} or \S 7.1 of \cite{CiLa11}. Let $\beta$ be the solution to the following equation, which involves a Legendre function: $P_{1/2}(\cos \pi/\beta) = 0$. Its value $\beta \simeq 1.3771$, and we have
\begin{prop}
\label{dimXY}
The singular spaces $\Xvec^S_{(k)}(\omega)$ and $\Yvec^S_{(k)}(\omega)$ are of finite dimension, namely
\begin{itemize}
\item For $k=0$
\begin{eqnarray*}
&&\dim \Yvec^S_{(k)}(\omega) := N_B = \mbox{ number of reentrant edges, }\\
&&\dim \Xvec^S_{(k)}(\omega) := N_E = N_B+ \mbox{ number of conical points with vertex angle } > \frac{\pi}{\beta}\,. 
\end{eqnarray*}
\item For $k \neq 0$
\begin{eqnarray*}
&& \dim \Yvec^S_{(k)}(\omega) := N_B = \dim \Xvec^S_{(k)}(\omega) := N_E = \mbox{ number of reentrant edges. }
\end{eqnarray*}
\end{itemize}
\end{prop}

\noindent By introducing now $(\xvec^{k}_{S,j})_{1 \leq j \leq N_E}$ and $(\yvec^{k}_{S,j})_{1 \leq j \leq N_B}$ the basis of $\Xvec^S_{(k)}(\omega)$ and $\Yvec^S_{(k)}(\omega)$ for a given Fourier mode $k$, we get that the singular parts of the Maxwell's equations solution can be decomposed into
$$
\Evec_S^{k}=\sum_{j=1}^{N_E} k^j_E \,\xvec^{k}_{S,j}\quad \mbox{ and } \quad
\Bvec_S^{k}=\sum_{j=1}^{N_B} k^j_B\, \xvec^{k}_{S,j} \,,
$$
where $k^j_E$ and $k^j_B$ are constants we will have to determine. This will be detailed in Section \ref{NumericResults}.\\

\noindent We present now the characterization of the singular basis $(\xvec^{k}_{S,j})_{1 \leq j \leq N_E}$ and $(\yvec^{k}_{S,j})_{1 \leq j \leq N_B}$. For simplicity, in what follows,  we will assume that there is only one singularity, that is $N_E=N_B=1$, and we will drop the index $j$. In these conditions, we are looking for the equations satisfied by $\xvec^{k}_{S}$ and $\yvec^{k}_{S}$. Following  \cite{AsCL02} \S5.1 and \S5.2, \cite{AsCLS03} \S3.1 or \cite{CiLa11} \S7.1, we obtain that they can be characterized {\em via} their variational formulation. Indeed, they are solution to the following homogeneous formulations
\begin{itemize}
\item For the space $\Xvec^S_{(k)}(\omega)$, the basis  $\xvec^k_{S}\in  \Xvec_{(k)}^{S}(\omega)$ solves:
\begin{equation}
\label{xsVF}
\left\{\begin{matrix}
a_{k}  \left ( \xvec^k_{S},\Fvec \right )=0,\; \; \forall \Fvec\in  \Xvec_{(k)}^{R}(\omega)& \\ 
\xvec^k_{S}\times \nvec|_{\gamma_b} = 0,\;  & \\
\xvec^k_{S}\cdot \nvec|_{\gamma_a} = 0,\;  & 
\end{matrix}\right.
\end{equation}

\item For the space $\Yvec^S_{(k)}(\omega)$, the basis  $\yvec^k_{S}\in  \Yvec_{(k)}^{S}(\omega)$ solves:
\begin{equation}
\label{ysVF}
\left\{\begin{matrix}
a_{k}  \left ( \yvec^k_{S},\Cvec \right )=0,\; \; \forall \Cvec\in  \Yvec_{(k)}^{R}(\omega)& \\ 
\yvec^k_{S}\cdot \nvec|_{\gamma} =0,\; \; \gamma:=\gamma_{b}\cup \gamma_{a} & 
\end{matrix}\right.
\end{equation}
\end{itemize}

\noindent As a consequence of Prop.\ref{propstab}, we readily get that the spaces $\Xvec^S_{(k)}(\omega)$ and $\Yvec^S_{(k)}(\omega)$ are satisfying 
$$
\Xvec^S_{(k)}(\omega)=\Xvec^S_{(2)}(\omega), \qquad \Yvec^S_{(k)}(\omega)=\Yvec^S_{(2)}(\omega), \quad \mbox{ for } |k| \geq  2\,.
$$
As recalled above, this stabilization property has a fundamental consequence on the numerical method, based on the decomposition of $\Xvec_{(k)}(\omega)$ and $\Yvec_{(k)}(\omega)$ in a regular and singular subspace: it will be sufficient to compute the singular basis $(\xvec_{S,j}^k)_{j=1,N_S}$ and $(\yvec_{S,j}^k)_{j=1,N_S}$ only for $|k| \leq 2$ and not for all $k$, the modes $\pm 2$ also serving as a non-orthogonal complement for the modes $|k| > 2$. More details and numerical illustrations will be given in Section \ref{NumericResults}. Another choice would be to derive, for each $k$, a mode-specific orthogonal basis. The interested reader will find a comparison for the Poisson equation in \cite{KiKw09}.

\section{Computation of singular basis}\label{ComputSingu}
\subsection {The case of $\xvec^k_{S}$}\label{subsecxsk}

\noindent We now present the numerical method to compute the singular part $\xvec^k_{S}$. The basis $\xvec^k_{S}\in  \Xvec_{(k)}^{S}(\omega)$ solves the singular and homogeneous problem (\ref{xsVF}) in $\omega$. Consequently,  if we try to solve it with a standard finite element method, we will get a zero solution.\\

\begin{figure}[htbp!]
\begin{tabular}{lr}
{
\includegraphics[scale=0.5]{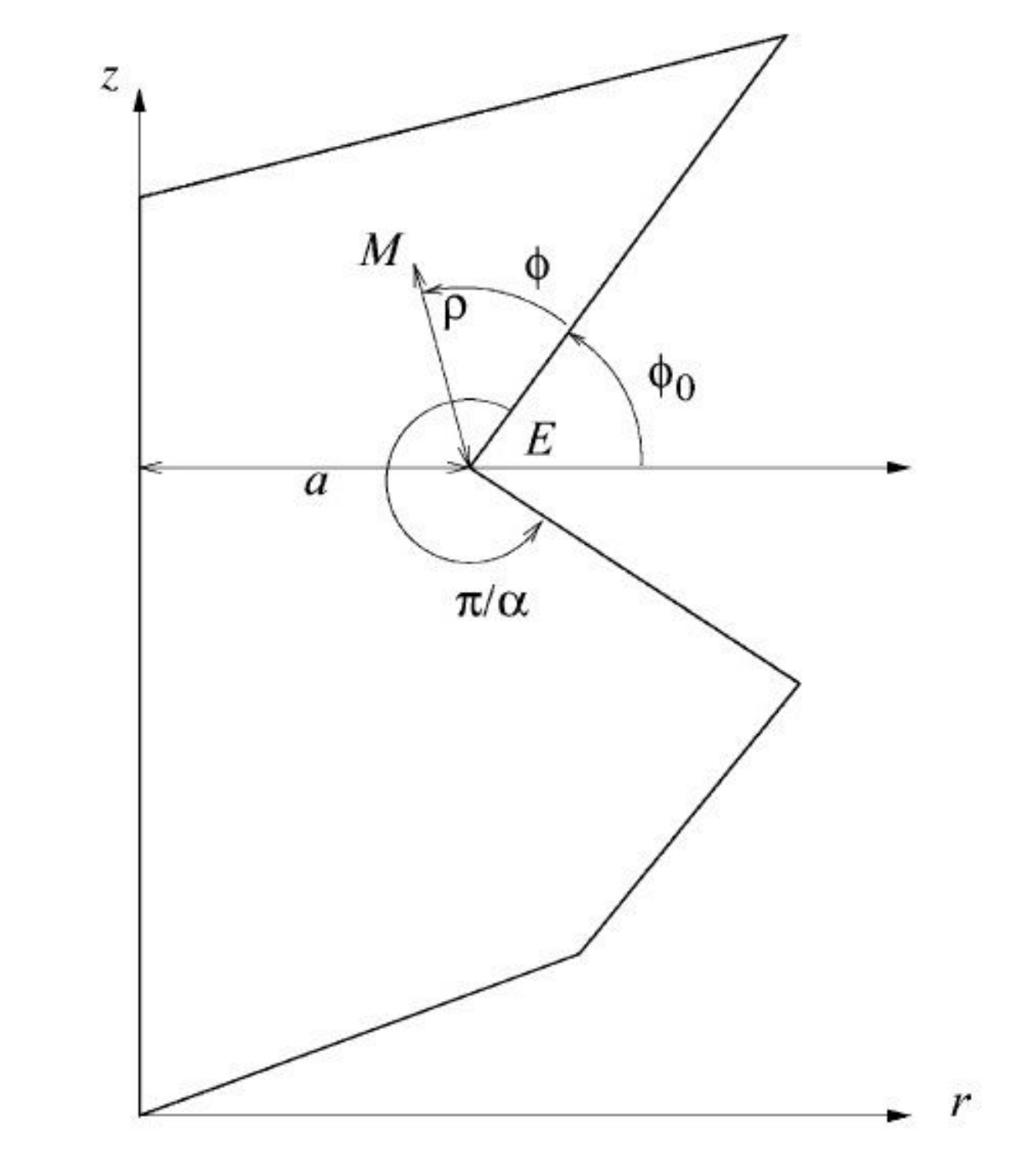}
}
&
\hspace*{2.cm}
{
\includegraphics[scale=0.5]{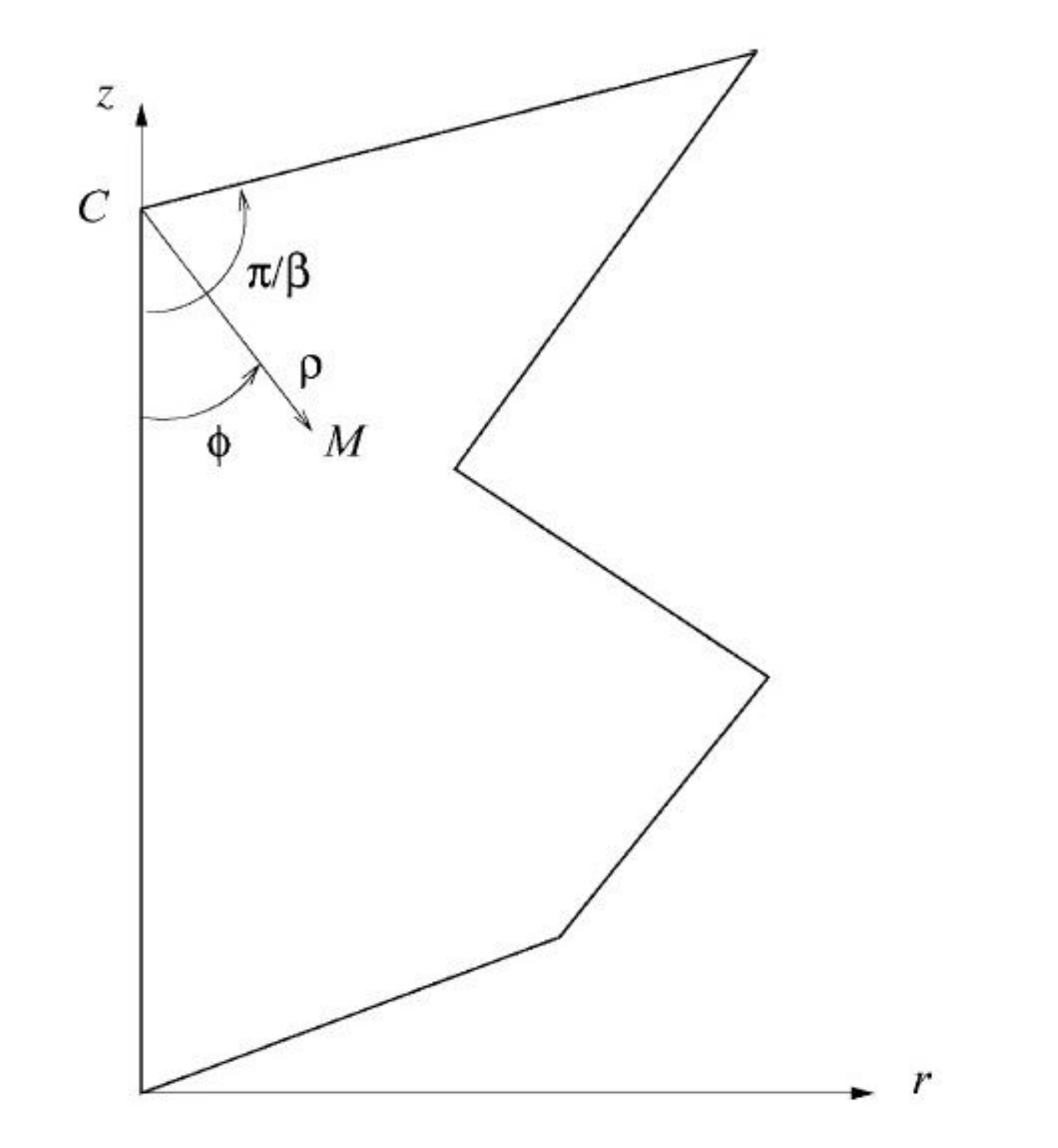}
}
\end{tabular}
 \caption{Local coordinates near a reentrant edge (left) and a conical point (right) .}
  \label{fig:local_coord}
\end{figure}

\noindent  To overcome this difficulty,  we rather use that the so-called {\em principal part }of the singularity $\Svec$ - the part that makes $\xvec^k_{S}$ singular - does not depend on the Fourier mode $k$. In these conditions, $\xvec^k_{S}$ can be decomposed into
$$
\xvec^k_{S}=\xvec^k_{S,reg}+\Svec
$$
where $\xvec^k_{S,reg} \in \Xvec^R_{(k)}(\omega)$ denotes the regular part of $\xvec^k_{S}$, that can be computed by a classical finite element method. Note also that in the electric case, as recalled in Prop.\ref{dimXY}, there exist two kinds of geometrical singularities $\Svec$:
\begin{enumerate}
\item for all $k \in \Z$, the one that exists in the neighborhood of a reentrant edge of $\Omega$, that is a reentrant corner of $\omega$, ($E$ in Figure \ref{fig:local_coord}-left), that we will denote $\Svec_e$,
\item only for $k=0$, the conical singularity ($C$ in Figure \ref{fig:local_coord}-right), that we will denote $\Svec_c$, that exists in the neighborhood of a conical vertex with an aperture greater than the limit vertex angle $\pi/\beta$ for $\beta \simeq 1.3771$ (so that $\pi/\beta \gtrsim 130^{\circ}43'$). 
\end{enumerate}

\noindent Figure \ref{fig:local_coord} shows the notations associated to these singularities. In particular, $(\rho,\phi)$ denotes the local polar coordinates centered at the reentrant edge $E$, the corresponding angle being called $\pi/\alpha$, $1/2 < \alpha < 1$. For the conical point $C$, $(\rho,\phi)$ are the local polar coordinates centered at this point, with the origin of $\phi$ on the $z$-axis.\\

\noindent  In these conditions,  the principal part $\Svec_e$ can be written as 
$\Svec_e=-\left ( r/a \right )\gradv_{0}\left [ \rho^{\alpha}\sin\left ( \alpha \phi \right ) \right ]$, 
whereas the principal part at the conical point (if any) can be expressed  as 
$\Svec_c=-\gradv_{0}\left [ \rho^{\nu}P_{\nu}\left ( \cos \phi \right ) \right ]$.
As above, $P_{\nu}$ denotes also here the Legendre function of index $\nu $, where $\nu \in ]0,1/2[$ is the index such that
$P_{\nu}(\cos(\pi/\beta))=0$.
The expressions of the  $\Svec_e$ and $\Svec_c$ in the basis $(\evec_r, \evec_\theta,\evec_z)$ are given by
\begin{equation}
\label{expressionofS_e}
\Svec_e=-\frac{r}{a}\alpha\rho^{\alpha-1} \begin{pmatrix}
\sin((\alpha-1) \phi-\phi_0)\\ 
0\\ 
\cos((\alpha-1) \phi-\phi_0)
\end{pmatrix}\,,
\end{equation}
\begin{equation}
\label{expressionofS_c}
\Svec_c=\nu\rho^{\nu-1} \begin{pmatrix}
 P_{\nu}\left ( \cos\phi\right )\cos\phi-P^1_{\nu}\left ( \cos\phi\right )\sin\phi\\ 
0\\ 
 P_{\nu}\left ( \cos\phi\right )\sin\phi +P^1_{\nu}\left (\cos\phi\right )\cos\phi
\end{pmatrix}\,,
\end{equation}

\noindent Remark that the term $r/a$ in (\ref{expressionofS_e}) has been introduced  to impose the boundary condition on the axis $r=0$, and can be viewed as a particular cut-off function. Moreover, the singularity $\Svec_c$ appearing only for $k=0$, it is related to a ``fully  axisymmetric" case, which was already treated in \cite{AsCLS03}, and will not be considered in the following (included numerical examples of Section \ref{NumericResults}).\\

\noindent Under these circumstances, we can compute  $\xvec^k_{S,reg}$ by solving the following variational formulation
$$
 \left\{\begin{matrix}
a_{k}  \left(\xvec^k_{S,reg},\vvec \right )=-a_{k}  \left ( \Svec_e,\vvec \right )\; \; \forall \vvec\in \Xvec_{(k)}^{R}(\omega), \\ 
\xvec^k_{S,reg}\times \nvec|_{\gamma_b} =-\Svec_e\times \nvec|_{\gamma_b},\;  &  \\
\xvec^k_{S,reg}\cdot \nvec|_{\gamma_a} =0,\;
\end{matrix}\right.
$$

\noindent The right-hand side of this equation can be computed analytically, by using the following
expressions of $\curlv_k \Svec_e$ and $\div_k \Svec_e$, involved in $a_{k}  \left ( \Svec_e,\vvec \right )$, see (\ref{akexpression}):
$$
\begin{matrix}
\curlv_k \Svec_e=\ds\frac{\alpha}{a}\,\rho^{\alpha-1} \begin{pmatrix}
-ik\cos( (\alpha-1 )\phi -\phi_0)\\ 
\cos( (\alpha-1 )\phi-\phi_0)\\ 
ik\sin(( \alpha-1)\phi - \phi_0)
\end{pmatrix} ,& 
\div_k \Svec_e=-\ds\frac{2\alpha}{a}\,\rho^{\alpha-1}\sin((\alpha-1)\phi-\phi_0)\,.
\end{matrix}
$$

\begin{Rem}
\noindent By construction, see for instance Eqs.(\ref{xsVF})-(\ref{ysVF}), the singular complement is orthogonal with respect to the form $a_k(\cdot,\cdot)$. Now, for the numerical implementation, it is also possible to orthonormalise the basis $\xvec^k_{S,j}$, for each singularity $j$, and to compute the basis vectors $(\xvec^{k,\perp}_{S,j})_j$ which are orthogonal to one another and to the regular space ${\textbf{X}}^R_{(k)}(\omega)$ (for $|k| \leq 2$). The same is true for the magnetic case. \\
Now, concerning the principal part of the singularity $\Svec$, it is the same, whatever the singular basis that we consider, the orthogonalisation process modifying only the regular part of $(\xvec^k_{S,j})_j$, and not $\Svec$. Computing such an orthogonal complement requires an additional computational effort on the one hand. On the other hand, the variational formulations (\ref{akexpressionE})-(\ref{akexpressionB}) are easier to solve because they contains fewer terms, some canceling due to orthogonality.
\end{Rem}

\subsection {The case of $\yvec^k_{S}$}\label{subsecysk}
\noindent Let us  turn our attention to the computational method for the singular part $\yvec^k_{S}$. This time, the basis  $\yvec^k_{S}\in  \Yvec_{(k)}^{S}(\omega)$ solves the system of equations (\ref{ysVF}).
This problem is singular and homogeneous in $\omega$: if we try to solve it with a standard finite element approach, we will get a zero solution.\\

\noindent To overcome this difficulty,  we rather use, as for the case of $\xvec^k_{S}$, that the principal part of the singularity $\Svec$, i.e. the part that makes $\yvec^k_{S}$ singular, does not depend on the Fourier mode $k$. In these conditions, $\yvec^k_{S}$ can be decomposed into
$$
\yvec^k_{S}=\yvec^k_{S,reg}+\Svec
$$
where $\yvec^k_{S,reg}$ denotes the regular part of $\yvec^k_{S}$, that can be computed by a classical finite element method. The expression of $\Svec$ in the basis $(\evec_r, \evec_\theta,\evec_z)$ is given by
\begin{equation}
\label{expressionofS}
\Svec=-\frac{r}{a}\alpha\rho^{\alpha-1} \begin{pmatrix}
\cos ((\alpha-1)\phi -\phi_0)\\ 
0\\ 
-\sin((\alpha-1)\phi -\phi_0)
\end{pmatrix}\,,
\end{equation}
where  the term $r/a$ in (\ref{expressionofS}) is useful to impose the boundary condition on the axis $r=0$, and is, here again, a particular cut-off function. Note also that in the magnetic case, there is no singularity due to the presence of conical vertex.\\

\noindent In these conditions, the function $\yvec^k_{S,reg}$ will be computed by solving the following variational formulation
\begin{equation}\label{FVyskstat}
 \left\{\begin{matrix}
a_{k}  \left(\yvec^k_{S,reg},\vvec \right )=-a_{k}  \left ( \Svec,\vvec \right )\; \; \forall \vvec\in \Yvec_{(k)}^{R}(\omega), \\ 
\yvec^k_{S,reg}\cdot\nvec|_{\gamma} =-\Svec \cdot \nvec|_{\gamma}\; \; .& 
\end{matrix}\right. 
\end{equation}

\noindent The right-hand side of this equation is computed analytically by using the following expressions of $\curlv_k \Svec$ and $\div_k \Svec$, involved in $a_{k}  \left ( \Svec,\vvec \right )$, see (\ref{akexpression}):
$$
\begin{matrix}
\curlv_k \Svec=\ds\frac{\alpha}{a}\,\rho^{\alpha-1} \begin{pmatrix}
ik\sin\ ((\alpha-1)\phi -\phi_0)\\ 
-\sin ((\alpha-1)\phi -\phi_0)\\ 
ik\cos ((\alpha-1)\phi -\phi_0)
\end{pmatrix} ,& 
\div_k \Svec=-\ds\frac{2\alpha}{a}\,\rho^{\alpha-1}\cos((\alpha-1)\phi -\phi_0)\,.
\end{matrix}
$$
\vspace*{0.2cm}

\noindent For the practical purpose of the computation, it is useful to express the bilinear form $a_k\left (\uvec,\vvec \right )$, depending on the values of $k$. Indeed, recall that, in our approach, we will have to compute the singular basis $\xvec^k_{S}$ and $\yvec^k_{S}$ only for $|k| \leq 2$. Performing a simple integration by parts shows that
\begin{eqnarray}\label{akuv}
a_k\left (\uvec,\vvec \right )&=& a_{0}\left ( \uvec_m,\vvec_m \right )+k^{2}\left ( \frac{\uvec_m}{r},\frac{\vvec_m}{r} \right )
+ \left ( \curlv u_\theta,\curlv v_\theta \right )+k^{2}\left ( \frac{u_\theta}{r},\frac{v_\theta}{r} \right ) \nonumber \\
&+&\imath k B\left(\uvec,\vvec \right )\,\,\,\,+\imath k C\left(\uvec,\vvec \right )\,,
\end{eqnarray}
where $a_{0}(\cdot,\cdot)$ denotes the operator $a_{k}(\cdot,\cdot)$ for $k=0$ (namely in the ``fully"  axisymmetric case), $\uvec_m:=\left ( u_r,u_z \right )$, the vector $\curlv$ of a scalar field $w$ being defined by  
$$
\curlv w:=-\partial_{z}w\evec_r+r^{-1}\partial_{r}\left ( rw \right )\evec_z\,.
$$
In addition, the two bilinear forms $B\left(\uvec,\vvec \right )$ and $C\left(\uvec,\vvec \right )$ are defined by 
$$
B\left(\uvec,\vvec \right ):=\int_{\gamma_b}\left (\uvec_m\cdot \nvec  \right )\bar{v}_\theta-u_\theta \left (\bar{\vvec}_m\cdot \nvec  \right )d\gamma\,,
$$
and 
$$
C\left(\uvec,\vvec \right ):=\int \int_\omega 2 \left(u_\theta\bar{v}_r  - u_r \bar{v}_\theta \right ) \frac{d\omega}{r}\,.
$$
Remark first that the term $B\left(\uvec,\vvec \right )$ vanishes as soon $\uvec\times\nvec = \vvec\times\nvec=0$ or $\uvec\cdot\nvec = \vvec\cdot\nvec=0$, that is exactly the case for the electric or magnetic fields (and also for $\xvec^k_{S}$ and $\yvec^k_{S}$), due to the perfect conductor boundary condition. Note also that the term $C\left(\uvec,\vvec \right )$ is not singular despite the presence of $1/r$ in the integral. Indeed, it is assumed that components $u_r, u_\theta$ and $v_r, v_\theta$ all belong to the space such that $\ds\int \int_\omega \frac{u v}{r} d\omega < \infty$, which ensures the convergence of the integral defining $C\left(\uvec,\vvec \right )$. However,  the numerical evaluation of this integral deserves special care, as detailed in \cite{AsCLS03}, \S 4.3.
\section{Numerical results} \label{NumericResults}

\noindent In this section, we present numerical results of singular field computations in axisymmetric domain. For the sake of simplicity, we will consider a 3-D top hat domain $\Omega$ with a reentrant circular edge, that corresponds, for a given $\theta$,  to an L-shaped 2-D domain $\omega$ with only one singular point, i.e. a reentrant corner (see figures below). To compute the numerical solutions, we introduce an unstructured mesh of $\omega$ made up of triangles, with no particular mesh refinement near the corner.  We then approximate the variational formulations presented in the sections above by using a $\mathrm{P}_1$ finite element method with FreeFem++ package,  which
implements a finite element method in space \cite{Hech12}.\\

\subsection {Computation of the singular basis for $|k| \leq 2$}
In what follows, we present numerical results for the magnetic case described in Subsection \ref{subsecysk}. The electric case (Subsection \ref{subsecxsk}) can be dealt in a similar way.
\subsubsection{Mode $k=\pm 1$}
\begin{figure}[htbp!]
\centering
\begin{tabular}{llr}
\includegraphics[scale=0.3]{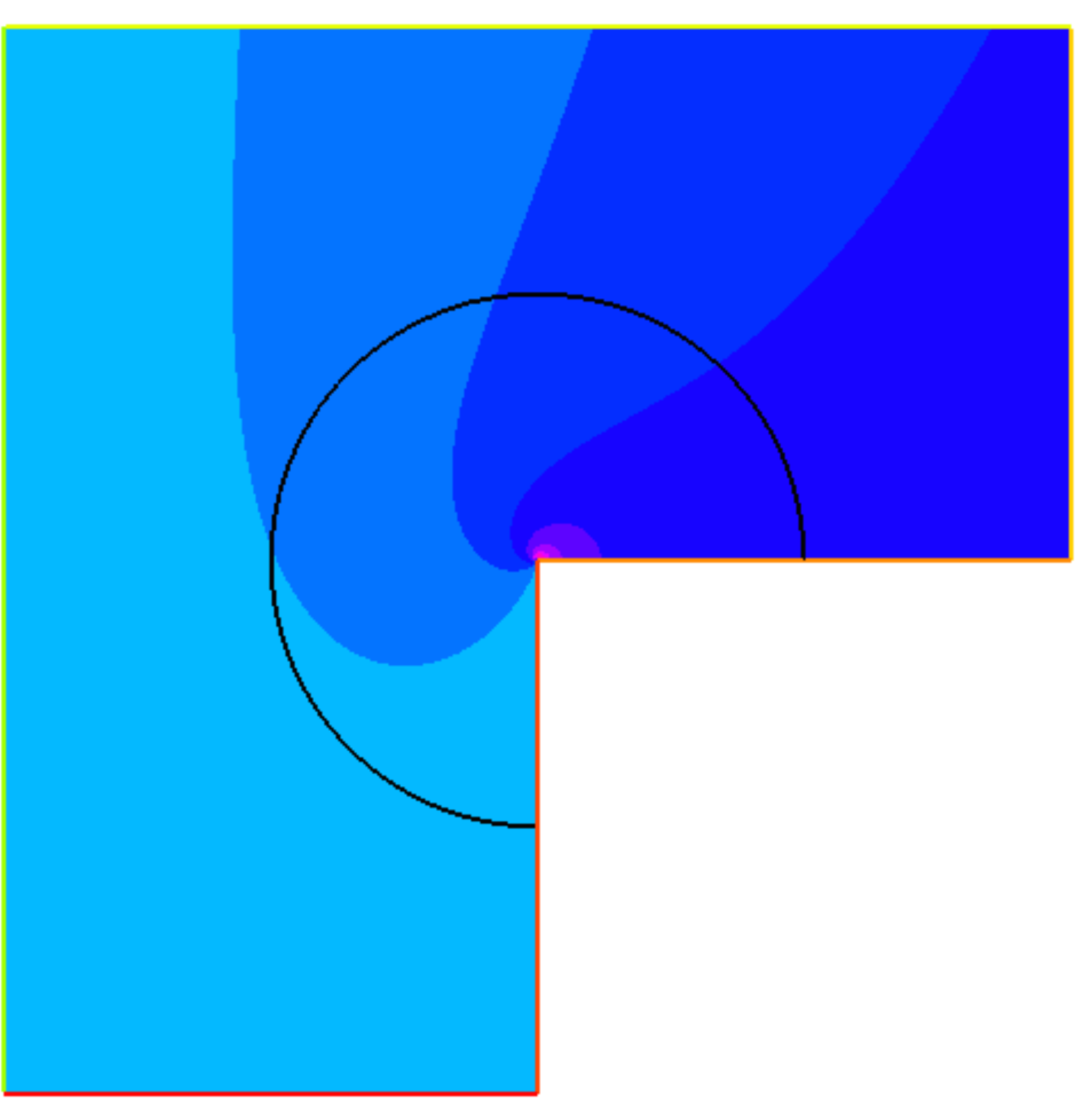} 
&
\hspace*{1.cm}
&
\includegraphics[scale=0.3]{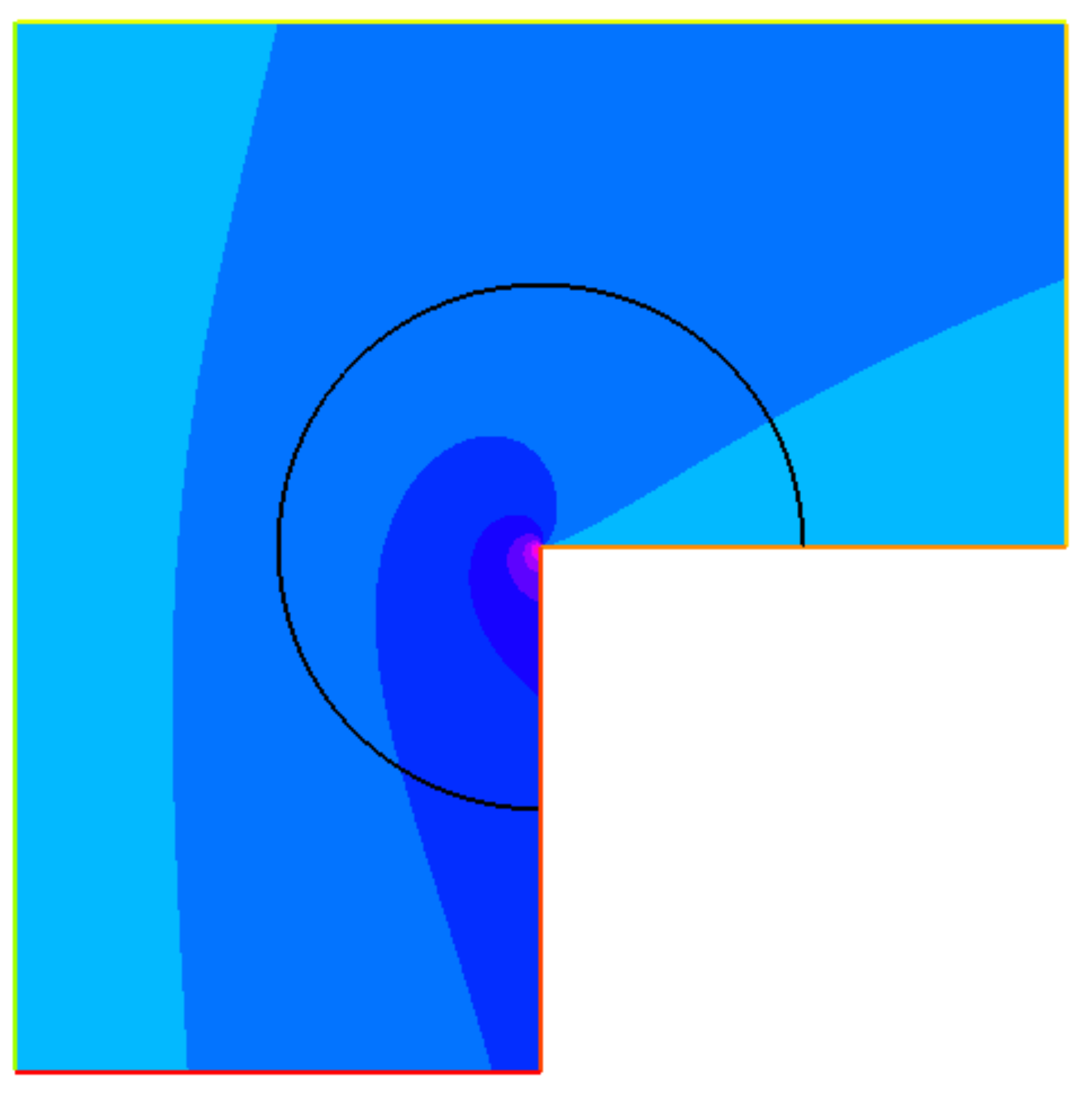}
\end{tabular}
 \caption{$S_r$ and $S_z$ components of the principal part of the singularity $\Svec$.}
  \label{fig:resultsA}
\end{figure}

We begin by presenting numerical results obtained by computing the singular basis $\yvec_S^k$ for $k=1$ and $k=-1$ respectively. For this purpose, we follow the procedure presented in Subsection \ref{subsecysk}. Below are depicted (see Figure \ref{fig:resultsA}) the non-vanishing components of the principal part of the singularity $\Svec$ introduced in (\ref{expressionofS}).

\begin{figure}[htbp!]
\centering
\begin{tabular}{llr}
\includegraphics[scale=0.3]{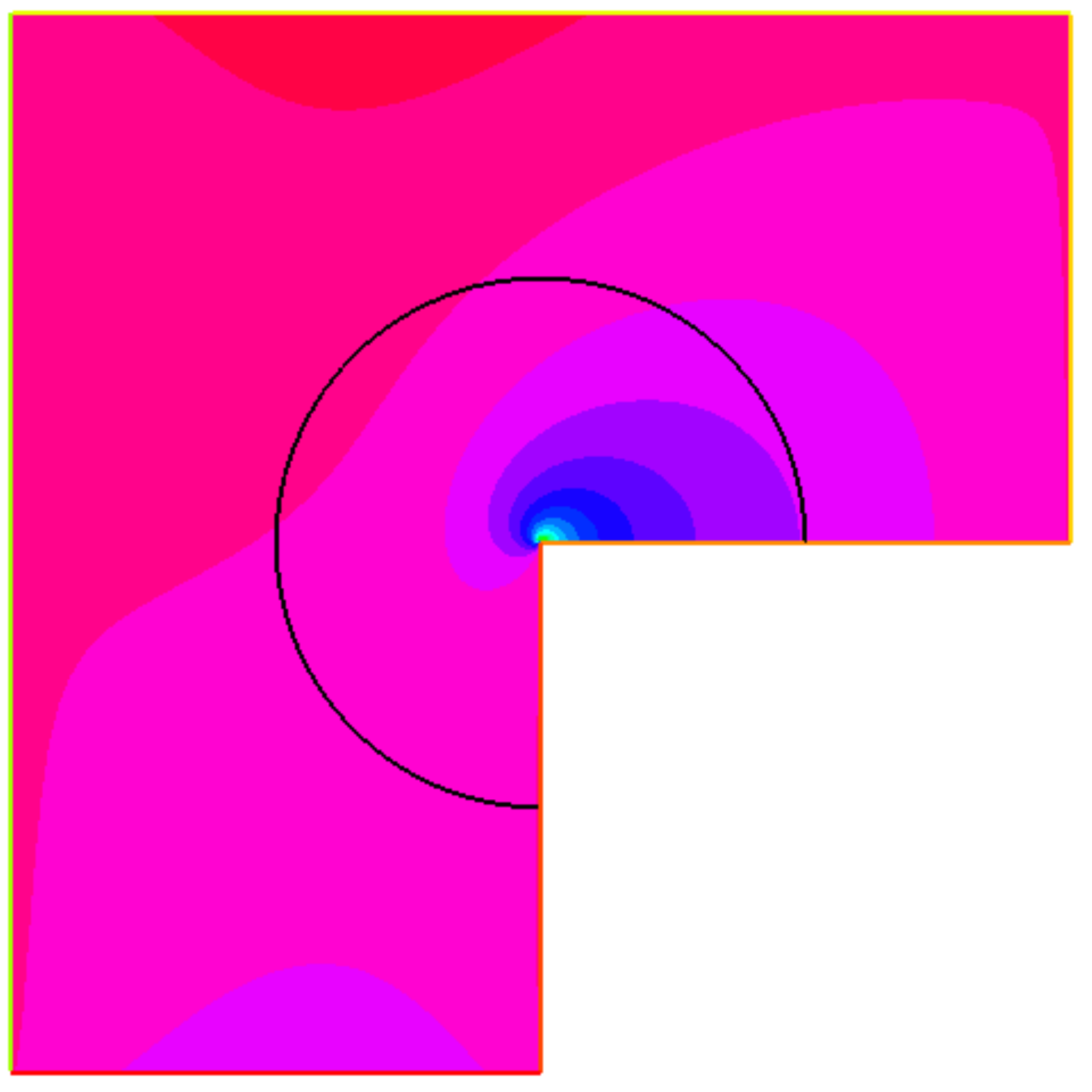}
&
\hspace*{1.cm}
&
\includegraphics[scale=0.3]{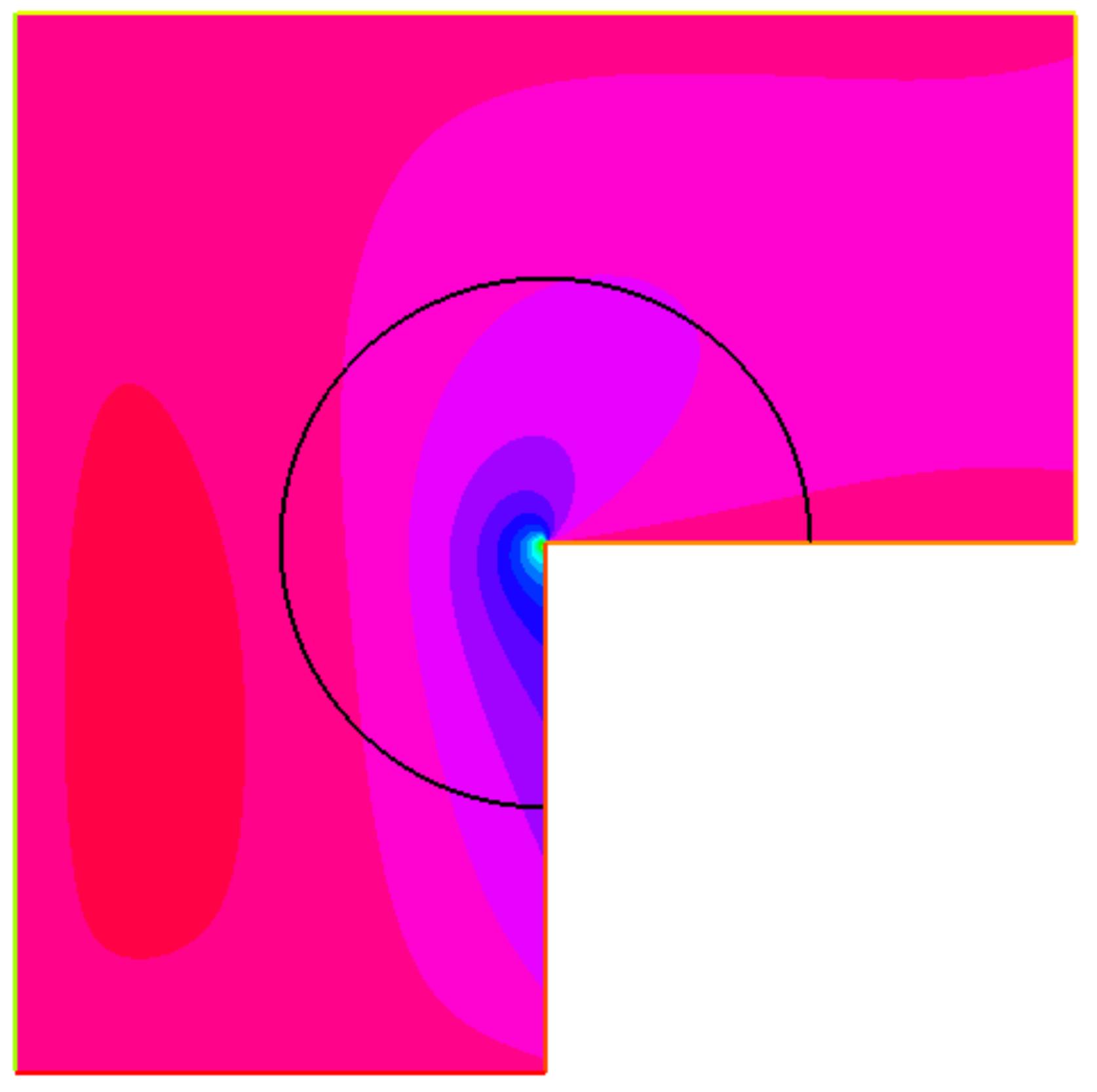}
\end{tabular}
 \caption{Real part $\Re(y_{S,r}^1)$ and $\Re(y_{S,z}^1)$ of the singular basis for $k=1$, $r$ and $z$ components.}
  \label{fig:resultsB}
\end{figure}

\begin{figure}[htbp!]
\centering
\begin{tabular}{llr}
\includegraphics[scale=0.3]{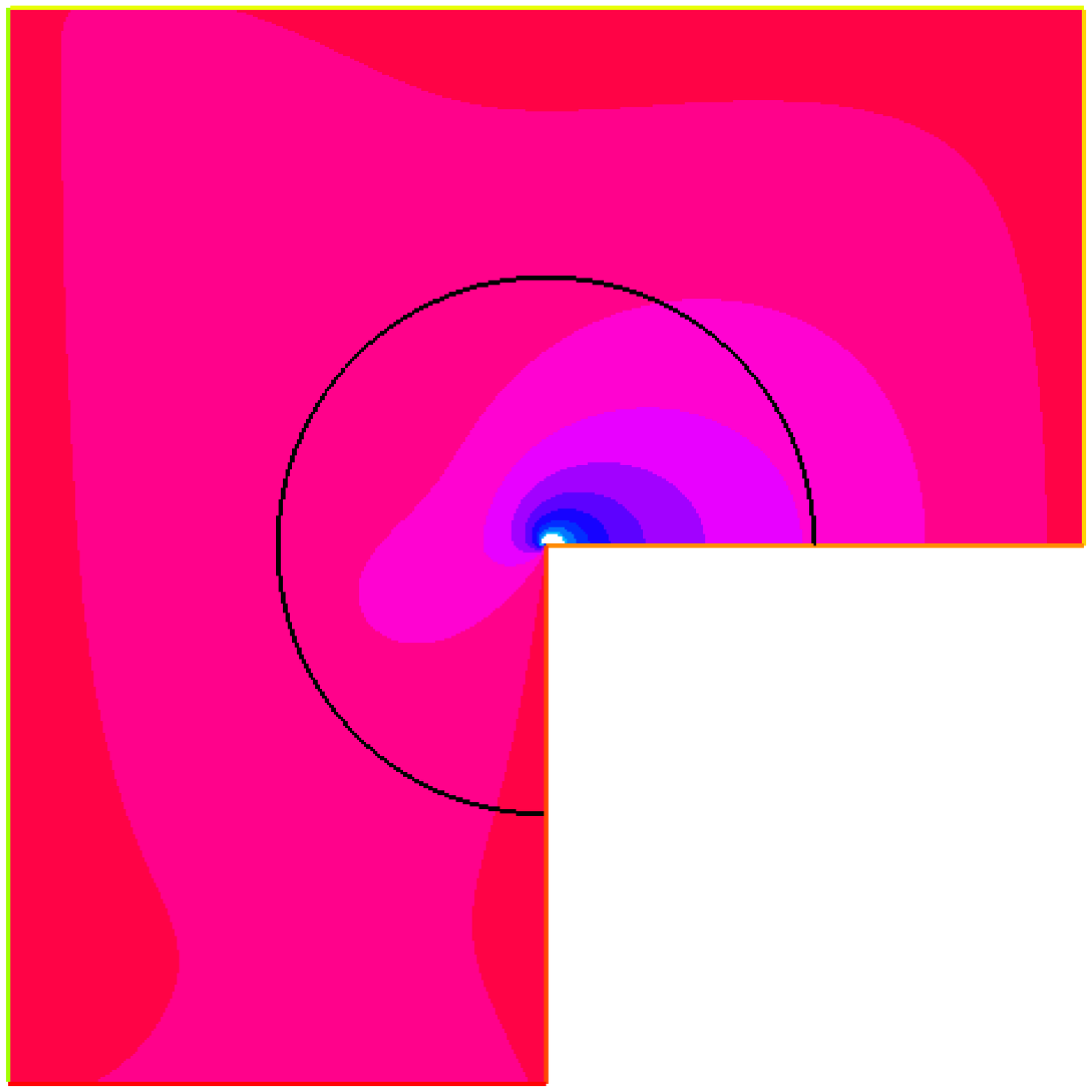}
&
\hspace*{1.cm}
&
\includegraphics[scale=0.3]{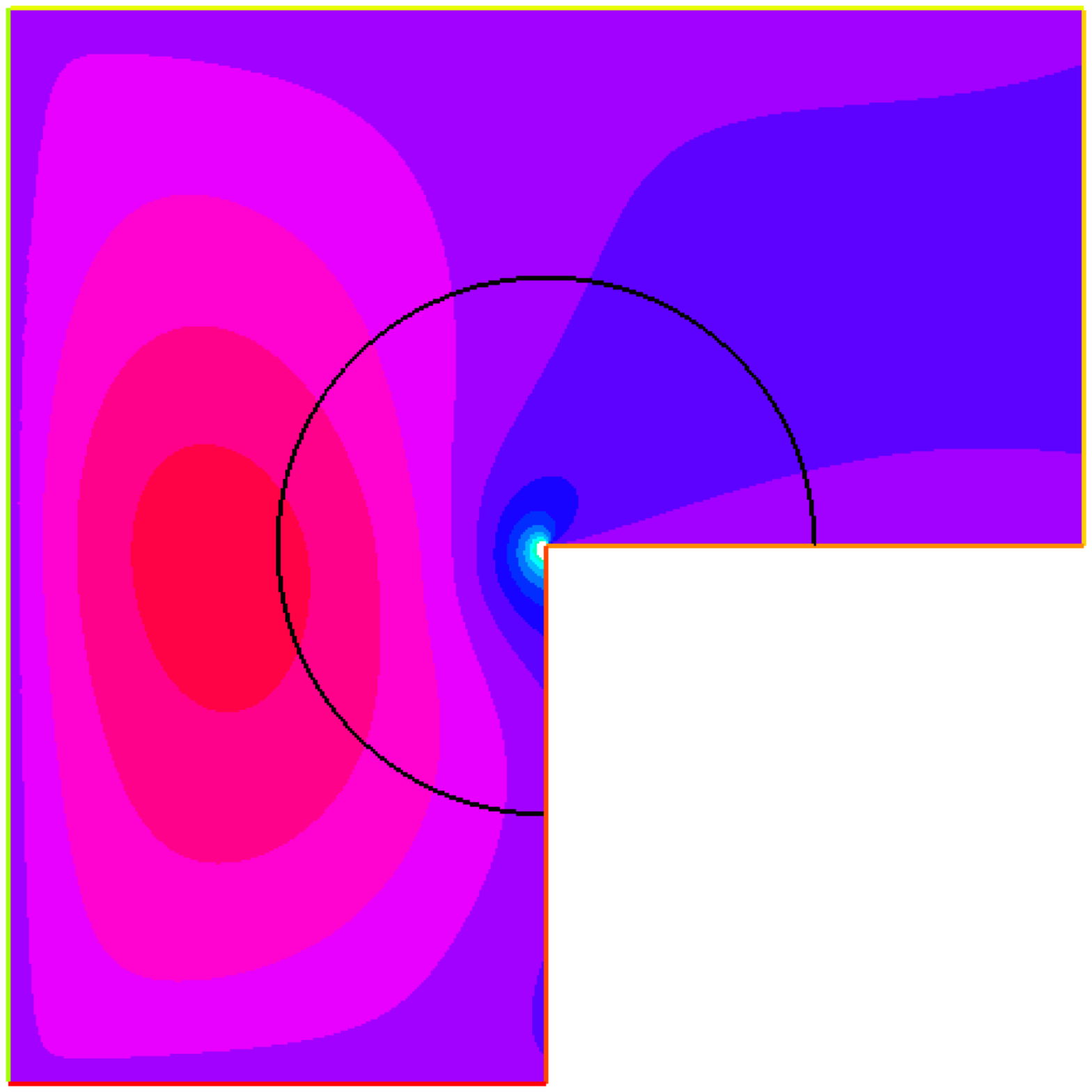}
\end{tabular}
 \caption{Real part $\Re(y_{S,r}^{-1})$ and $\Re(y_{S,z}^{-1})$ of the singular basis for $k=-1$, $r$ and $z$ components.}
  \label{fig:resultsBbis}
\end{figure}

\noindent We then compute the singular basis $\yvec_{S}^1$ and $\yvec_{S}^{-1}$ for $k=1$ and for $k=-1$. In Figures \ref{fig:resultsB} and \ref{fig:resultsBbis}, we depicted their real part (remember that they are complex quantities, as soon as $k \neq 0$) obtained by solving the variational formulation (\ref{FVyskstat}), with a standard $P_1$ finite element method. As one can see, the method is able to capture the singular behavior of the solution near the reentrant corner of $\omega$ (edge in $\Omega$), whereas a conforming finite element method can not yield such a result.
\subsubsection{Mode $k = \pm 2$}
\begin{figure}[htbp!]
\centering
\begin{tabular}{llr}
\includegraphics[scale=0.3]{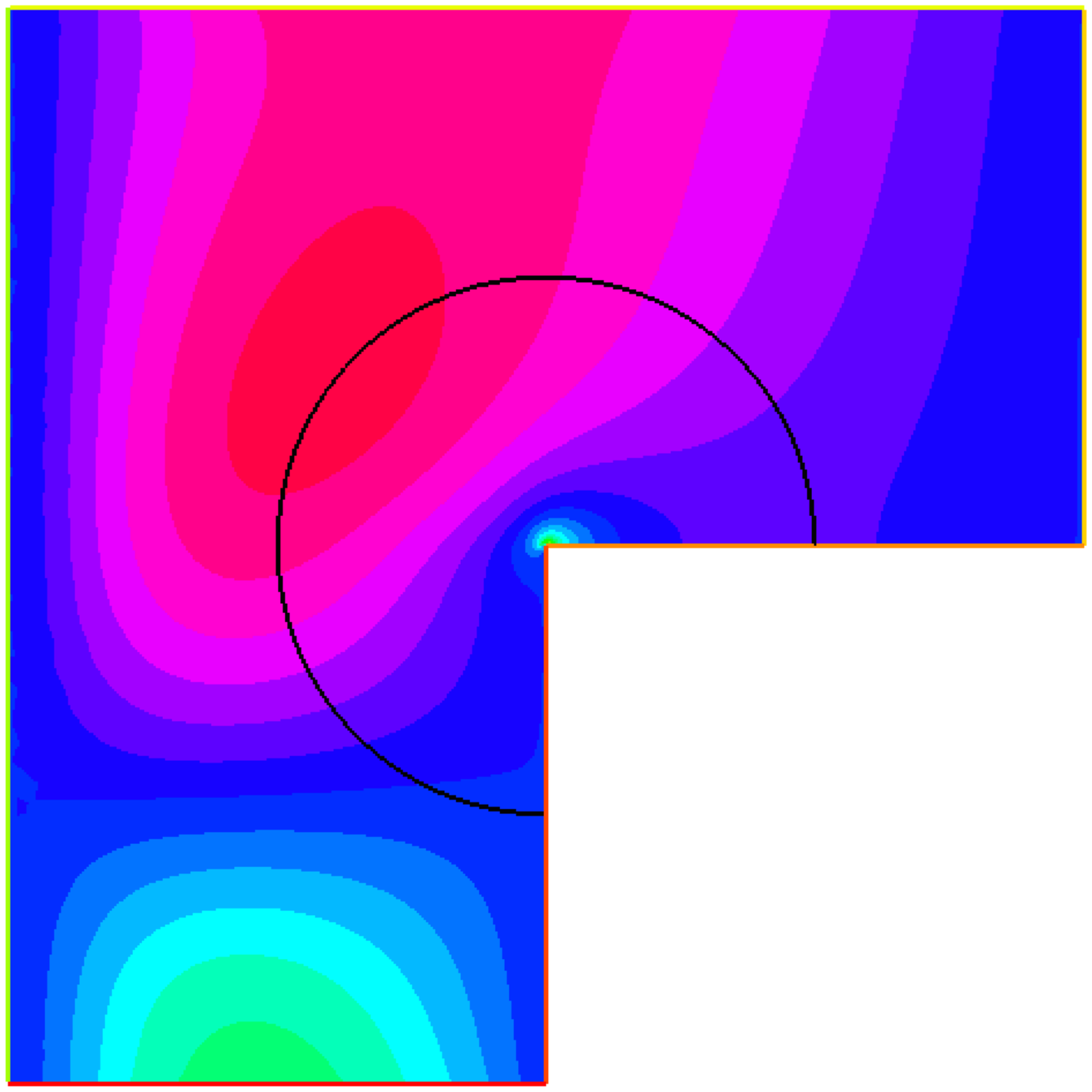}
\end{tabular}
 \caption{Real part $\Re(y_{S,r}^2)$  
 of the singular basis for $k=2$, $r$ component.}
  \label{fig:resultsC}
\end{figure}

For the sake of completeness, we present below the results obtained for $k=2$, see Figure \ref{fig:resultsC}. As explained above, there is no need in our numerical strategy to compute modes for $|k| \geq 2$, since the modes corresponding to $|k| = 2$ appear as fundamental modes for $|k| > 2$. In other words, they can be used to compute the solution for all $|k| \geq 2$, despite the loss of orthogonality. This point will be illustrated in the next subsection. The numerical method is exactly the same as for $k=\pm1$, and practically, we have only to change  the parameter $k$ in the formulation.
\subsubsection{Mode $k=0$}
\noindent In this subsection, we deal with the particular mode $k=0$, that corresponds to the ``fully" axisymmetric problem. In that case, we already proposed other ways to derive a method that can capture the singular solution, see for instance \cite{AsCLS03},  \cite{AsRa15}. It is interesting here to compare the results obtained by the two approaches.\\

\noindent  Let us first briefly recall the principle derived for the ``fully" axisymmetric case. As a first step, we look for $P_S$,  a non-vanishing, singular, harmonic function  solution to a homogeneous Laplace problem set in $\omega$, with a homogeneous Dirichlet boundary condition on the boundary of the domain. Note that $P_S$ is not equal to zero, since we are looking for a solution that has not enough  regularity to be variational. As a consequence, computing it with a standard variational formulation would give $P_S=0$ as a solution.\\

\noindent  In the same spirit as in Section \ref{ComputSingu}, we decompose $P_S$ in a singular principal part, explicitly known, and a regular part, that can be computed with a finite element method. Then, using isomorphisms proved in \cite{AsCL03}, we introduce the unique potential $\psi$ solution to the Laplace problem with $P_S$ as the right-hand side. Finally, the basis $\yvec_S^0$ can be computed by taking the $\curlv_0$ of the potential $\psi$.\\
\begin{figure}[htbp!]
\centering
\begin{tabular}{llr}
\includegraphics[scale=0.118]{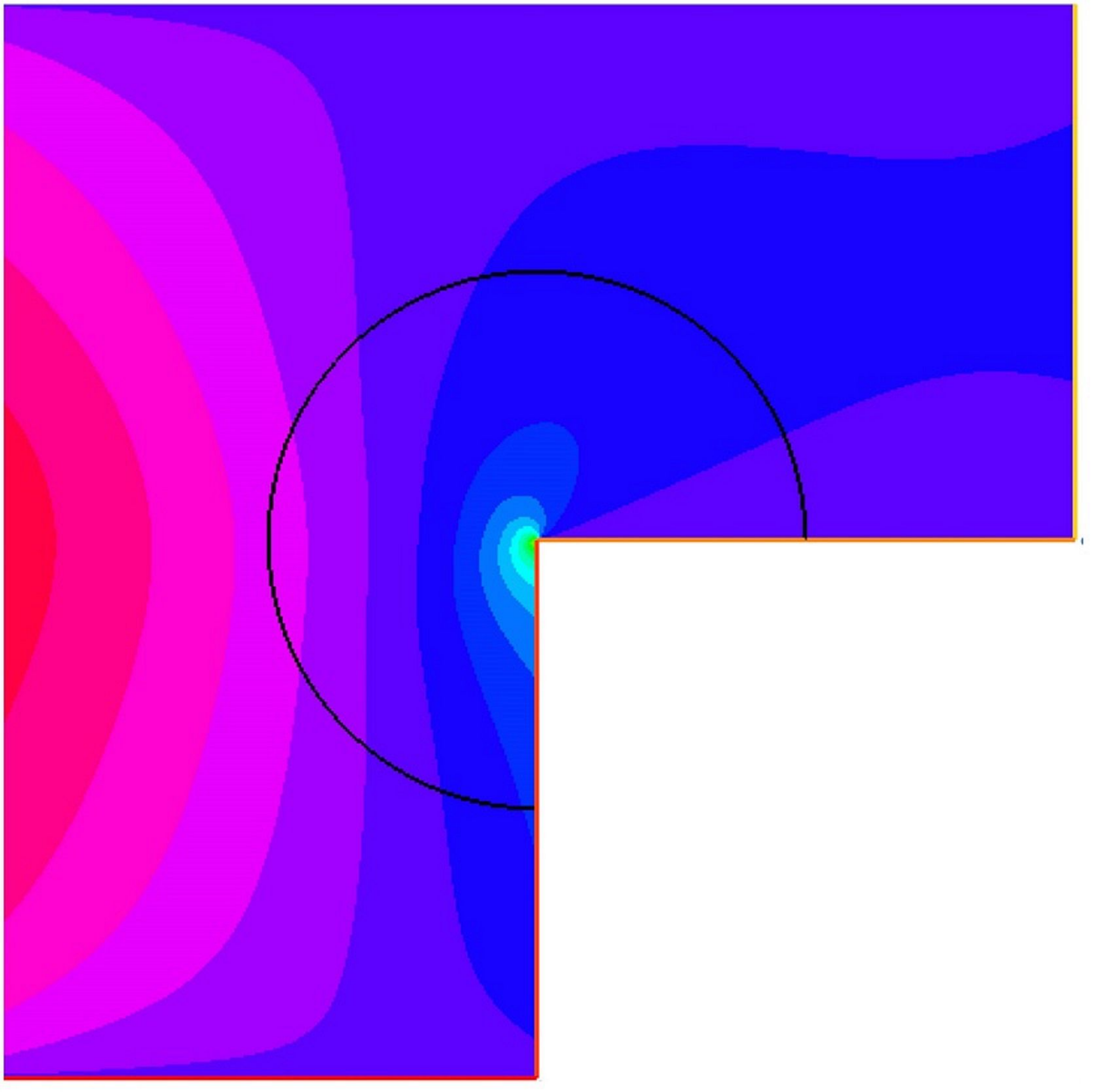}
&
&
\includegraphics[scale=0.3]{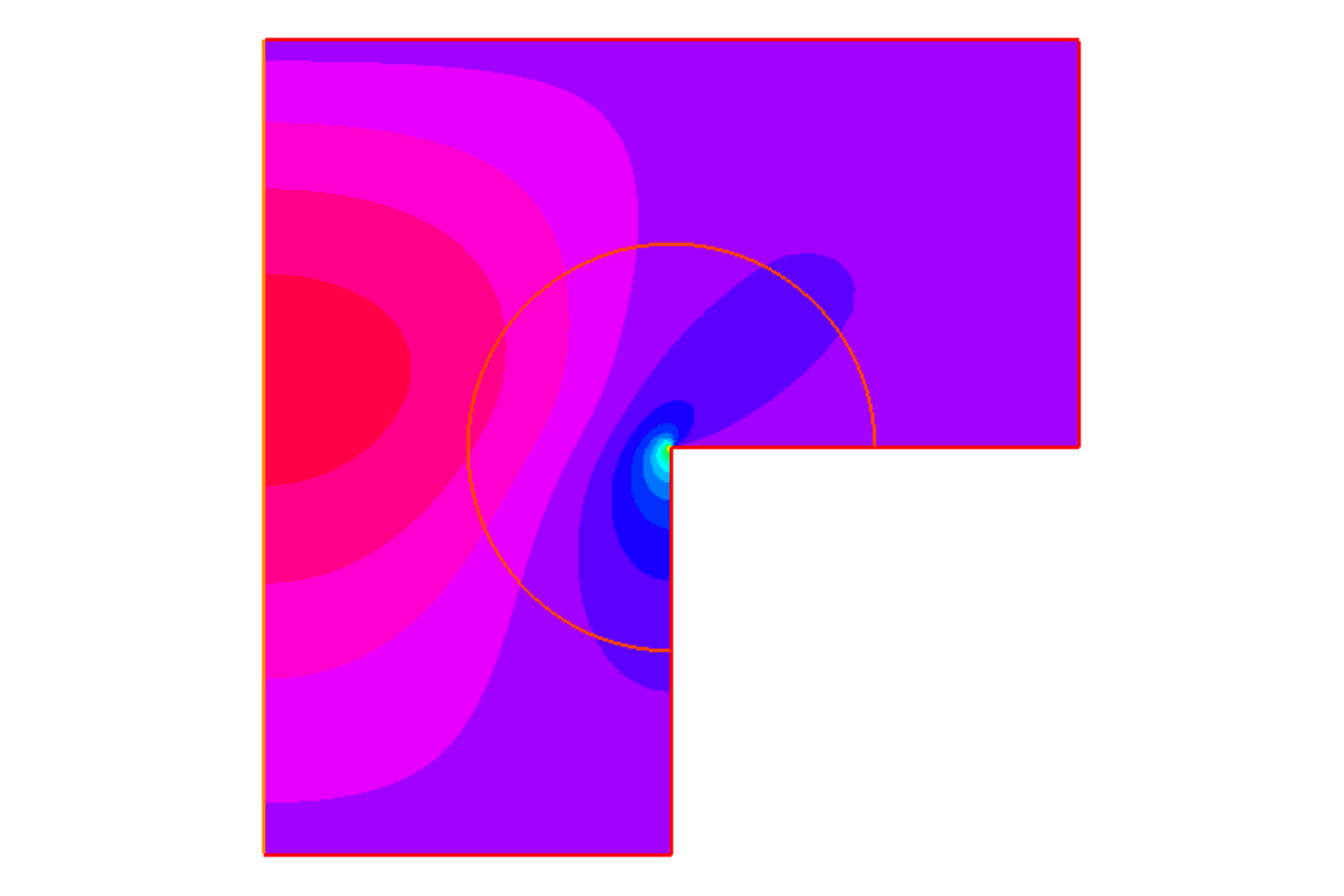}
\end{tabular}
 \caption{ $\yvec^0_{S}$ (left) compared with the ``fully"  axisymmetric  solution obtained in \cite{AsRa15} (right), $z$ component.}
  \label{finalresult}
\end{figure}

\noindent The numerical results (from \cite{AsRa15}) are depicted in Figure \ref{finalresult}-right and have to be compared to the one obtained with the method proposed in this article, by choosing $k=0$. As one can see (the colors of the scale are not exactly the same), there is a good agreement between the results obtained by the two methods.\\

\begin{Rem}
\noindent Once the basis $\yvec_S^k$ have been obtained, one can compute the magnetic field $\Bvec^k $ by solving the system (\ref{akexpressionB}). For this purpose, one has first to use the direct decomposition (\ref{decompEBk}) that is written here
$$
\Bvec^k=\Bvec_R^k + \Bvec_S^k\,,
$$
together with the characterization of the singular part $\Bvec_S^k=C^k \yvec_S^k$, since, for each $k$, $\Bvec_S^k$ belongs to the one-dimensional subspace $\mbox{span}\{ \yvec_S^k \}$. We readily get
$$
\Bvec^k=\Bvec_R^k \oplus C^k \yvec_S^k\,,
$$
where $C^k $ is a real constant to be determined. To that aim, one can derive the variational formulation (\ref{akexpressionB}) in which the test function $\Cvec$ is chosen equal to $\yvec_S^k$, namely 
$$
a_{k}\left ( \Bvec_R^k + C^k \yvec_S^k, \yvec_S^k   \right )= \left (\mathbf{f}^k_\Bvec,{\mathbf{curl}}_k\,  \yvec_S^k \right )+  \left (g^k_\Bvec,\mathrm{div}_k \yvec_S^k \right ) \,.
$$
Then,  using that $\mathrm{div}_k \Bvec^k=g^k_\Bvec=0$, $\mathrm{div}_k  \yvec_S^k$ is known (here is equal to $0$) together with orthogonality between  ${\mathbf{curl}}_k \,\Bvec_R^k$ and ${\mathbf{curl}}_k \, \yvec^{k}_{S}$, for a given $|k| \leq 2$, we have
\begin{equation}
\label{expr-Ck}
C^k=\dfrac{\left (\mathbf{f}^k_\Bvec,{\mathbf{curl}}_k \,  \yvec_S^k \right )}{\left ({\mathbf{curl}}_k \,  \yvec^{k}_{S}, {\mathbf{curl}}_k \,  \yvec^{k}_{S} \right ) }\,.
\end{equation}

\noindent Hence, the singular field $\Bvec_S^k$ being known, one can compute the regular part of the magnetic field $\Bvec_R^k$ still using the variational formulation with a regular test function  belonging to $\Yvec^R_{(k)}$. With the same arguments as above, one directly  obtains that the regular part $\Bvec_R^k$ solves the formulation
\begin{equation}\label{VFBregk}
a_{k}\left ( \Bvec_R^k,\Cvec  \right ):= \left ({\mathbf{curl}}_k  \Bvec_R^k , {\mathbf{curl}}_k  \Cvec  \right )+\ \left (\mathrm{div}_k \Bvec_R^k, \mathrm{div}_k \Cvec  \right )= \left (\mathbf{f}^k_\Bvec,{\mathbf{curl}}_k \, \Cvec \right ), \,\,\, \forall \Cvec \in \Yvec^R_{(k)}\,.
\end{equation}
Since all the field involved in (\ref{VFBregk}) are regular, this can be computed by a finite element method. The total field is then ``rebuilt" by addition, for each $|k| \leq 2$, then by using expansion (\ref{FourierTrunc}).
\end{Rem}

\noindent Following the previous remark, these computations can be easily adapt to compute the electric field $\Evec^k$.

\subsection {Computation of a Fourier mode for $|k| > 2$}

\noindent In the following, we compute a Fourier mode for  $|k| > 2$, here  $k = 3$, to assess the use of the singular basis of the mode $k=2$ for determining the modes for $|k| > 2$. Indeed, owing to the stabilization property (Prop. \ref{propstab}), one can use the singular basis $\yvec_S^2$ as a fundamental mode for all $|k| >  2$. This makes the method effective in the sense that it reduces the number of singular modes to be calculated.\\

\noindent  Nevertheless, $\yvec_S^2$ serving as a non-orthogonal complement for $|k| > 2$, some terms in the variational formulation no longer cancel each other out, and the method must take into account this loss of orthogonality. Another possibility could be to derive a ``mode-specific orthogonal" basis for each mode $k$. In this article, we derive the first approach, as proposed in \cite{CJKLZ05I} for the Poisson problem, the additional computational effort (i.e. the non-vanishing terms due to the loss of orthogonality) being small, and not significantly changing the implementation of the method.\\

\noindent Let us use again the magnetic case as an illustration. For $|k| > 2$, the total magnetic field $\Bvec^k=\Bvec_R^k + \Bvec_S^k=\Bvec_R^k + C^k \yvec_S^2$ is computed by solving a coupled system, the unknowns being $\Bvec_R^k$ and $C^k$. Indeed, using the same arguments as in Remark 2, one gets, taking successively $\Cvec \in \Yvec^R_{(k)}$ and  $\yvec_S^2$ $ \in \Yvec^S_{(2)}$ as test functions in (\ref{akexpressionB})
\begin{equation}
\label{systkgt2}
\left\{
 \begin{array}{l}
a_{k}\left ( \Bvec_R^k, \Cvec) + C^k a_{k}(\yvec_S^2, \Cvec \right )= \left (\fvec^k_\Bvec,\curlv_k  \Cvec  \right ) , \quad \forall \Cvec \in  \Yvec^R_{(k)}\,\\
\\
a_{k}\left ( \Bvec_R^k, \yvec_S^2) + C^k  a_{k}(\yvec_S^2,\yvec_S^2   \right )= \left (\fvec^k_\Bvec,\curlv_k \yvec_S^2 \right ) \,.
\end{array}
\right.
\end{equation}

\noindent After discretization in space, formulation (\ref{systkgt2}) can be expressed equivalently as a linear system
\begin{equation}
\label{systkgt2dis}
\left\{
\begin{array}{l}
  {\mathbb K^k_{rr}} {\Bvec^k_R} +  C^k \,  {\mathtt{Y}}_{RS}^2 = {\mathtt{F}}_{\Bvec,R}^k \;,\\
  \\
 ^t{\mathtt{Y}}_{RS}^2 \,{\Bvec^k_R} +  \alpha^2_k \, C^k  = f_{\Bvec,S}^2 \,.
\end{array}
\right.
\end{equation}
where ${\mathbb K^k_{rr}}$ is the matrix associated to the term $a_{k}\left (\Bvec_R^{k}, \Cvec  \right )$, ${\mathtt{Y}}_{RS}^2$ is a vector coming from the discretization of $a_{k}(\yvec_S^2, \Cvec )$ and $^t{\mathtt{Y}}_{RS}^2$ its transpose.  The term $\alpha^2_k$ is a scalar derived from the  discretization of $a_{k}(\yvec_S^2, \yvec_S^2)$, the vector ${\mathtt{F}}_{\Bvec,R}^k$ and the scalar $f_{\Bvec,S}^2$ coming from the right-hand sides of (\ref{systkgt2}).\\

\noindent  Hence, solving system (\ref{systkgt2dis}) amounts to computing ``together" the constant $C^k$ and the regular part $\Bvec^k_R$, that were solved separately for $|k| \leq 2$, where the $a_k$-orthogonality allows us to ``decouple" $C^k$ and $\Bvec^k_R$.  Hence,  the additional effort is not  very significant.\\

\begin{Rem}
\noindent It is worth noting that the knowledge of $a_2\left (\mathbf{u},\mathbf{v} \right )$, computed for $k=2$, can be used to compute the bilinear form $a_k(\mathbf{u}, \mathbf{v} )$ expressed in (\ref{akuv}), using the identity (for all $k$) 
$$
a_k \left (\mathbf{u},\mathbf{v} \right ) = a_2 \left (\mathbf{u},\mathbf{v} \right ) + (k^2-4)\,\left( \frac{\mathbf{u}}{r},\frac{\mathbf{v}}{r} \right ) +\imath (k-2) \,C \left(\mathbf{u},\mathbf{v} \right )
$$
\end{Rem}

\noindent To illustrate our approach and to assess the use of the singular basis of the mode $k=2$ for the other modes, we first compute 
for $k=3$ the magnetic total field $\Bvec^3:= \Bvec_R^3+ C^3\, \yvec_S^2$ by solving the system (\ref{systkgt2dis}). Then, we compute the singular basis $\yvec_S^3$, that allows one to compute differently $\Bvec^3$, this time with the method presented in remark 2, extended to $k=3$, i.e. using the $a_3$ orthogonality.\\

\begin{figure}[htbp!]
\centering
\begin{tabular}{llr}
\includegraphics[scale=0.3]{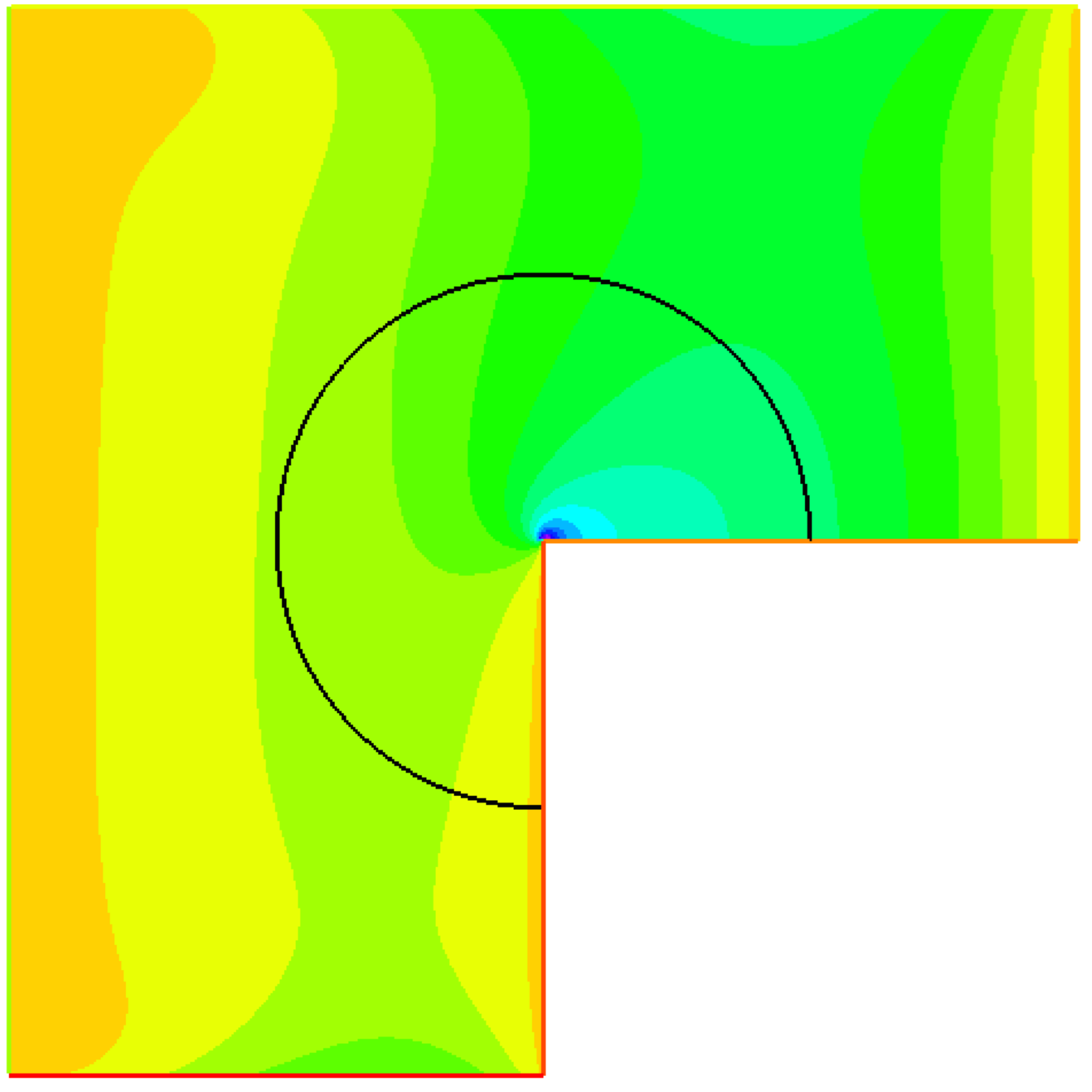}
&
\hspace*{1.cm}
&
\includegraphics[scale=0.3]{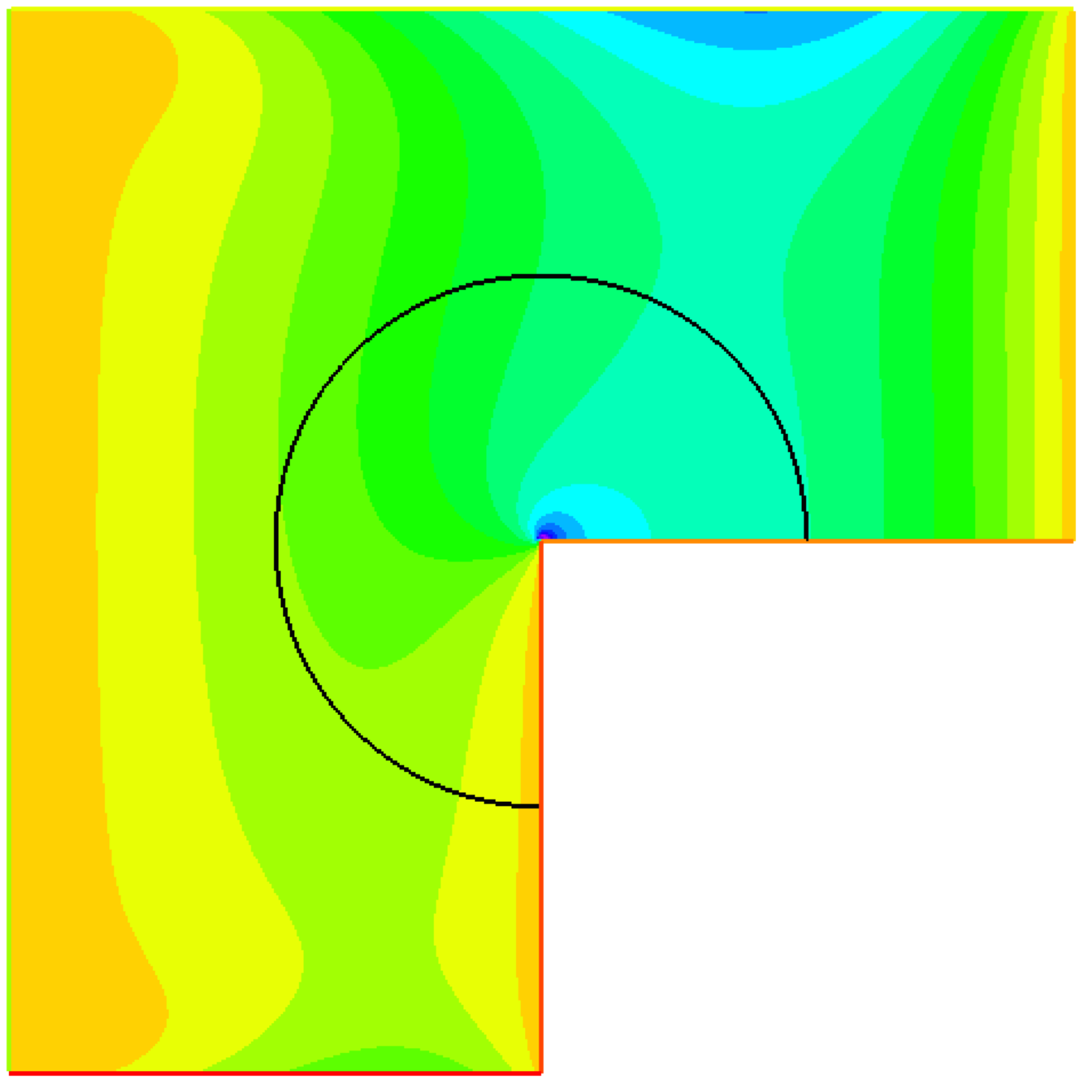}
\end{tabular}
 \caption{Real part $\Re(B_{r}^3)$ of the r-component of  $\Bvec^k$ for $k=3$, computed by two methods: (left) with $\yvec_S^3$ and the $a_3$ orthogonality,  (right) using the singular basis $\yvec_S^2$ as a ``fundamental" mode.}
  \label{fig:resultsk3}
\end{figure}

\noindent In Figure \ref{fig:resultsk3}, we compare the results obtained by the two approaches. The real part $\Re(B_{r}^3)$ (r-component) of  $\Bvec^3$ using the singular basis $\yvec_S^3$ compared to the method using $\yvec_S^2$ as a ``fundamental" mode. Note that both pictures are displayed with the same scale for better comparison. As one can see, the two methods, orthogonal and non-orthogonal,  that are equivalent on the continuous level, give very similar results,  but with slight differences, that can be explained in the following way.\\

\noindent When using the orthogonal approach, namely computing $\Bvec^k=\Bvec_R^k \oplus C^k \yvec_S^k$, one basically uses that terms like $\left ({\mathbf{curl}}_k \, {\Bvec^k_R}, {\mathbf{curl}}_k \,  \yvec^{k}_{S} \right )$ or $\div_k \yvec^{k}_{S}$ are exactly equal to 0, yielding expression (\ref{expr-Ck}) of $C^k$. Now, even on the discrete level, these terms are exactly equal to 0 and are not computed, since they are not involved in (\ref{expr-Ck}). The same is true when one computes numerically the regular part $\Bvec_R^k$ by solving (\ref{VFBregk}).\\

\noindent Now, in the non-orthogonal approach (\ref{systkgt2}) (or in (\ref{systkgt2dis}) after discretization), similar terms, that are different from 0, are actually computed so that the precision of the numerical approach (mesh, finite element used, etc.) leads to small differences in the numerical results. This can be improved simply by using a more refined mesh, or a more accurate finite element method (for instance a $\mathrm{P}_2$ finite element method instead of a $\mathrm{P}_1$ one) depending on the degree of accuracy one needs.\\

\noindent On the other hand, note that the orthogonal method is quite sensitive to the value of $C^k$, which is obtained after dividing by the singular term $\left ({\mathbf{curl}}_k \,  \yvec^{k}_{S}, {\mathbf{curl}}_k \,  \yvec^{k}_{S} \right )$. In that case, this denominator has to be numerically approach very carefully.
\section{Conclusion} \label{Conclu}
We aimed to solve the three-dimensional static Maxwell equations in a singular axisymmetric domain. We therefore presented a numerical method that can be viewed as an extension to 3D axisymmetric problems of the Singular Complement Method.
The first step was to reduce the dimension by using a Fourier transform in the azimuthal variable $\theta$. This gives us 
a series of 2D Maxwell's equations, depending on the Fourier variable $k$.  The second step was to deal with the 2D singularity for each $k$. For this purpose, we basically used (for each $k$) a splitting of the space of solutions in a regular subspace, which is equal to the entire space when the domain is smooth or convex, and a singular subspace. Due to a stabilization property, we only need to compute the singular part for a few values of $k$. Finally, numerical examples have been proposed to illustrate the method. Currently, we devote our attention to the extension to the time-dependent problem. Potential extensions and applications to the nonlinear case and optimal control problems can also be considered, as for instance in \cite{Yous13}.

\end{document}